\documentclass[12pt]{article}
\usepackage{graphicx}
\usepackage{color}
\usepackage{amssymb,amsmath}
\def\g{g}
\def\z{\zeta}
\def\B{\mathcal{B}}
\def\BSP{\B_{\mathrm{sp}}}
\def\DSP{\D_{\mathrm{sp}}}
\def\C{\mathcal{C}}
\def\D{\mathcal{R}}
\def\F{\mathcal{F}}

\def\G{\mathcal{G}}
\def\H{\mathcal{H}}
\def\K{\mathcal{K}}

\def\P{\mathcal{P}}
\def\PP{\mathrm{PP}}

\def\R{\mathcal{R}}
\def\S{\mathcal{S}}
\def\T{\mathcal{T}}
\def\N{\mathcal{N}}
\def\M{\mathcal{M}}
\def\bc{\mathrm{bc}}
\def\tc{\mathrm{tc}}
\def\Aut{\mathrm{Aut}}
\def\aut{\mathrm{aut}}
\newcommand{\Ba}{\B_{\mathrm{all}}}
\newcommand{\Da}{\D_{\mathrm{all}}}
\newcommand{\Fa}{\F_{\mathrm{all}}}

\def\bc{\mathrm{bc}}
\def\tc{\mathrm{tc}}

\newcommand{\SPG}{\G_{\mathrm{sp}}}
\newcommand{\nn}{\mbox{\rm\small 1} \hspace{-0,30em} 1}
\newtheorem{theorem} {Theorem}

\newtheorem{defn}[theorem]{Definition}
\newtheorem{thm}[theorem]{Theorem}
\newtheorem{corol}[theorem]{Corollary}

\newtheorem{propos}[theorem]{Proposition}

\def\proof{\noindent{\bf Proof}\/.\ }

\begin{document}
\title{ 
Two-connected graphs with prescribed three-connected components}
\author{Andrei Gagarin\thanks{Jodrey School of Computer Science, Acadia University, Wolfville, Nova Scotia, Canada, B4P 2R6.}, Gilbert Labelle\thanks{Laboratoire de Combinatoire et d'Informatique Math\'ematique (LaCIM), Universit\'e du Qu\'ebec \`a Montr\'eal (UQAM), Montr\'eal, Qu\'ebec, CANADA, H3C 3P8. With the partial support of NSERC (Canada).}, Pierre Leroux\footnotemark[2],\\ and\\ Timothy Walsh\footnotemark[2]
}
\maketitle
\vspace{-5mm}
\begin{abstract}
We adapt the classical  3-decomposition of any 2-connected graph  to the case of  simple graphs  (no loops or multiple edges). 
By analogy with the block-cutpoint tree of a connected graph, 
we deduce from this decomposition a bicolored tree $\tc(g)$ associated with any 2-connected graph $g$, whose white vertices are the \emph{3-components} of $g$ (3-connected components or polygons) and whose black vertices are bonds linking together these 3-components, arising from separating pairs of vertices of $g$.  
Two fundamental relationships on graphs and networks follow from this construction.
The first one is a dissymmetry theorem which leads to the expression of the class $\B=\B(\F)$ of 2-connected graphs, all of whose 3-connected components belong to a given class $\F$ of 3-connected graphs, in terms of various rootings of $\B$.  The second one is a functional equation which characterizes the corresponding class $\D=\D(\F)$ of two-pole networks all of whose 3-connected components are in $\F$. 
All the rootings of $\B$ are then expressed in terms of $\F$ and $\D$.  There follow corresponding identities for all the associated series, in particular the edge index series. Numerous enumerative consequences are discussed. 
\end{abstract}
\section{Introduction} 
A graph is assumed to be simple, that is, with no loops or multiple edges. 
A graph $G$ is called \textit{$k$-connected} if at least $k$ of its vertices and their incident edges must be deleted to disconnect it.  In general, a $k$-connected graph is assumed to have more than $k$ vertices. However, by convention, the complete graph $K_2$ will be considered as a $2$-connected graph.
A connected graph $G$ is \textit{planar} if there exists a $2$-cell embedding (i.e. each face is homeomorphic to an open disk) of $G$ on the sphere, with similar definitions for \textit{toroidal} and \textit{projective-planar} graphs. 

A \textit{species} is a class $\C$ of labelled combinatorial structures (for example, graphs or rooted trees) which is closed under isomorphism. Each $\C$-structure has an underlying set (for example, the vertex set of a graph), and isomorphisms are obtained by relabelling along bijections between the underlying sets. Unlabelled structures are defined as isomorphism classes of structures. 
We sometimes denote a species by the name of representatives of the isomorphism classes. For example, $K_n$ is used to denote the species of complete graphs on $n$ vertices. One advantage of species is that very often combinatorial identities can be expressed at this structural level, within the algebra of species and their operations (sum, product, substitution and a special substitution, of networks for edges, etc.).  There follow corresponding identities for the various generating series that are used for labelled and$\slash$or unlabelled enumeration. The reader is referred to the book \cite{BLL} for more details on species, their operations and associated series.

\begin{figure}[h] 
\begin{center} \includegraphics[width=4.0in]{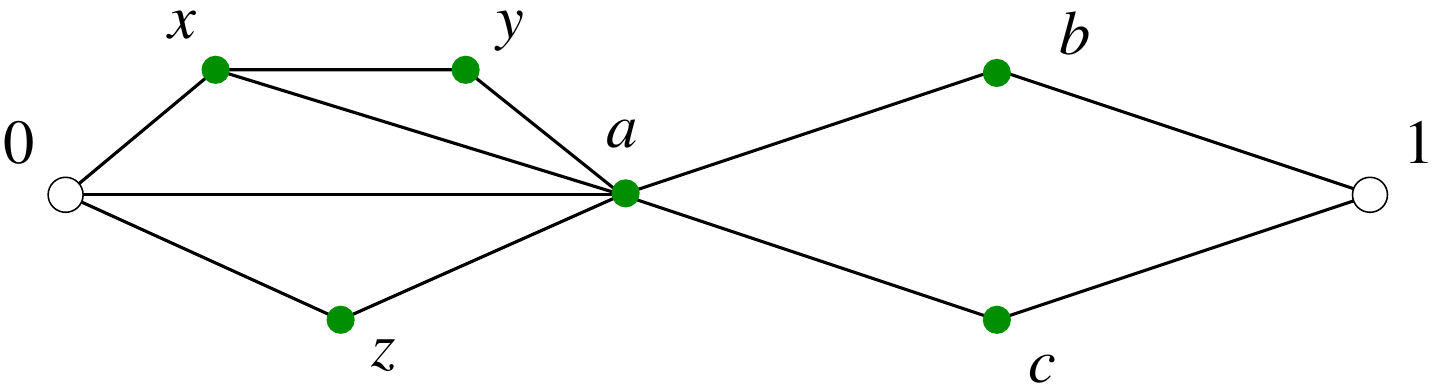}
\end{center}
\vspace{-4mm}
\caption{A 01-network   \label{fig:reseauxexemple}}
\end{figure}

The two-pole networks that we use have distinguished poles 0 and 1 and are called 01-networks. A 01-\textit{network} (or more simply a \emph{network}) is defined as a connected graph $N$ with two distinguished vertices $0$ and $1$, such that the graph $N\cup 01$ is $2$-connected, where the notation $N\cup ab$ is used for the graph obtained from $N$ by adding the edge $ab$ if it is absent. See Figure \ref{fig:reseauxexemple} for an example. The vertices $0$ and $1$ are called the {\it poles} of $N$, and all the other vertices of $N$ are called \textit{internal} vertices. The internal vertices of a network form its underlying set.  The trivial network, consisting of only the poles 0 and 1 and having no edges, is denoted by $\nn$.

In this paper, we first adapt to the case of simple graphs the classical 3-decomposition of 2-connected multigraphs (see Maclaine \cite{Maclaine} Tutte \cite{Tutte1,Tutte2}, Hopcroft and Tarjan \cite{Tarjan} and Cunningham and Edmonds \cite{Cunningham}).  By analogy with the block-cutpoint tree $\bc(c)$ of a connected graph $c$, we deduce from this decomposition a bicolored tree $\tc(g)$ associated with any 2-connected graph $g$, whose white vertices are the \emph{3-components} of $g$ (3-connected components or polygons)  and whose black vertices are bonds linking together these 3-components, arising from separating pairs of vertices of $g$, acting as hinges. 

Two fundamental relationships on graphs and networks follow from this construction. 
The first one is a dissymmetry theorem which leads to the expression of the class $\B=\B(\F)$ of 2-connected graphs, all of whose 3-connected components belong to a given class $\F$ of 3-connected graphs in terms of various rootings of $\B$. 
See Theorem \ref{theo:dissym} below. The second one is a functional equation which characterizes the corresponding class $\D=\D(\F)$ of non-trivial 01-networks all of whose 3-connected components are in $\F$. See Theorem \ref{theo:reseaux}.  Note that the 3-components of a network $N$ are obtained by considering the 3-decomposition of the graph $N\cup 01$.  

All the rootings of $\B$ are then expressible in terms of $\F$ and $\D$ and hence also $\B$ itself by virtue of the dissymmetry theorem. Although more or less implicit in previous work of one of the authors (see \cite{Timothy,Timothyunlabelled}), these identities are given here for the first time in the structural context of species.  There follow corresponding identities for all the associated series, in particular the edge index series, and numerous enumerative consequences are obtained. 
%

Among the examples that we have in mind and that will be discussed further in this paper are the following:
\begin{enumerate}
\item  If we take $\F=\Fa$, the class of all 3-connected graphs, then we have $\B(\F)=\Ba$, the class of all 2-connected graphs, and $\D(F)=\Da$, the class of all non-trivial 01-networks.
\item  If we take $\F=0$, the empty species, then $\B$ and $\D$ are the classes of 2-connected \emph{series-parallel} graphs and of \emph{series-parallel} networks, respectively.
\item
One of the motivations for the present paper was to extend some earlier tables for the number of $K_{3,3}$-free projective planar and toroidal 2-connected graphs (see \cite{GLL3}), which require the enumeration of \emph{strongly planar networks}, that is of non-trivial networks $N$ such that the graph $N\cup 01$ is planar. This class, denoted by $\N_P$, can be obtained by considering the class $\F=\F_P$ of planar 3-connected graphs.  Then the corresponding species $\B_P=\B(\F_P)$ is the class of planar 2-connected graphs, since a graph is planar if and only if all its 3-connected components are planar, and $\D(\F_P)=\N_P$, the class of strongly planar networks.
\item
As quoted by Thomas (see \cite{thomas}, Theorem 1.2, page 1), Wagner \cite{Wagner} has shown that a 2-connected graph is $K_{3,3}$-free if and only if it can be obtained from planar graphs and $K_5$'s by means of 2-sums (see also Kelmans \cite{Kelmans}). This means that if we take $\F = \F_P + K_5$, then the corresponding $\B=\B(\F)$ is the class of $K_{3,3}$-free 2-connected graphs. This fact was also observed by Gimenez, Noy and Ru\'e in \cite{Noy}.
\end{enumerate}
%

Section 2 contains the dissymmetry theorem.  Section 3 discusses various operations on 01-networks, in particular series and parallel composition and the substitution of networks for edges in graphs or networks. It also presents the fundamental relationship characterizing the species $\D=\D(\F)$ and the expressions of the various rootings of $\B$ in terms of $\F$ and $\D$. Applications to the labelled enumeration of these species are also presented.
Section 4 is devoted to the 
series techniques for species of graphs and networks that are necessary for their unlabelled enumeration.  These results are then applied to the enumeration of several classes of graphs and networks in Section 5. 
%
%
%
%
\section{A dissymmetry theorem for $2$-connected graphs.} 
In the standard decomposition of a 2-connected multigraph (multiple edges allowed but no loops) into 3-components (see \cite{Tutte1,Tutte2,Tarjan} and \cite{Cunningham}), the components are either 3-connected graphs (here called \emph{3-connected components}), polygons with at least three sides, or bonds, that is sets of at least 3 multiple edges.

\begin{figure}[h] 
\begin{center}
\includegraphics[width=4.4in]{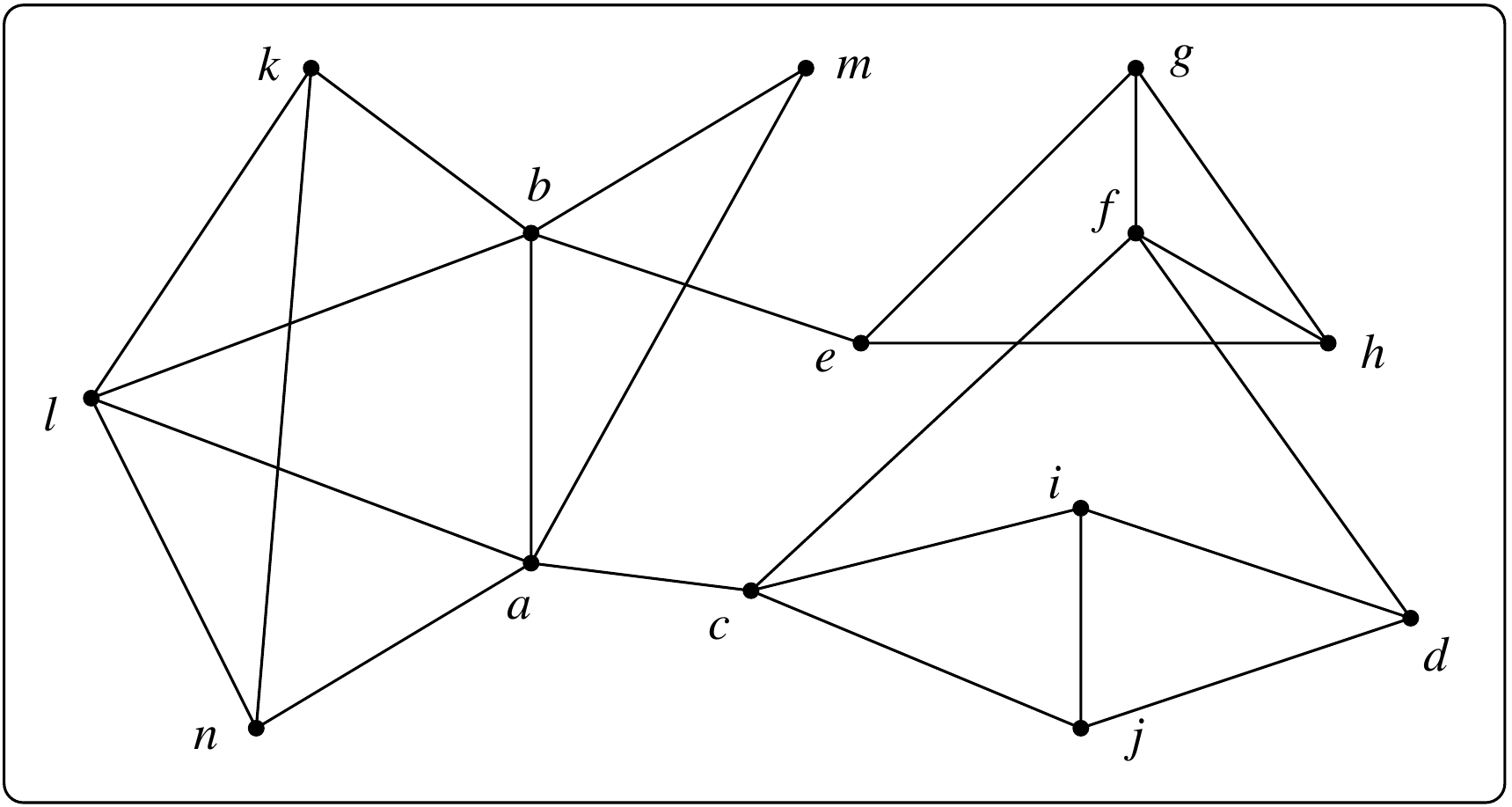}
\end{center}
\vspace{-4mm}
\caption{\small A 2-connected graph $g$
\label{fig:graphe2con}}
\end{figure}
%
However, in the case of simple graphs, bonds are not needed as 3-components and the decomposition is simpler.  We illustrate the construction with the graph $g$ of Figure \ref{fig:graphe2con}. By definition, a \emph{separating pair} of a 2-connected graph $g$ is a pair of vertices $\{x,y\}$ whose removal disconnects the graph. One can then in each resulting connected component re-introduce the two vertices $x$ and $y$ together with their  incident edges and the edge $xy$. 
\begin{figure}[h] 
\begin{center}
\includegraphics[width=4.6in]{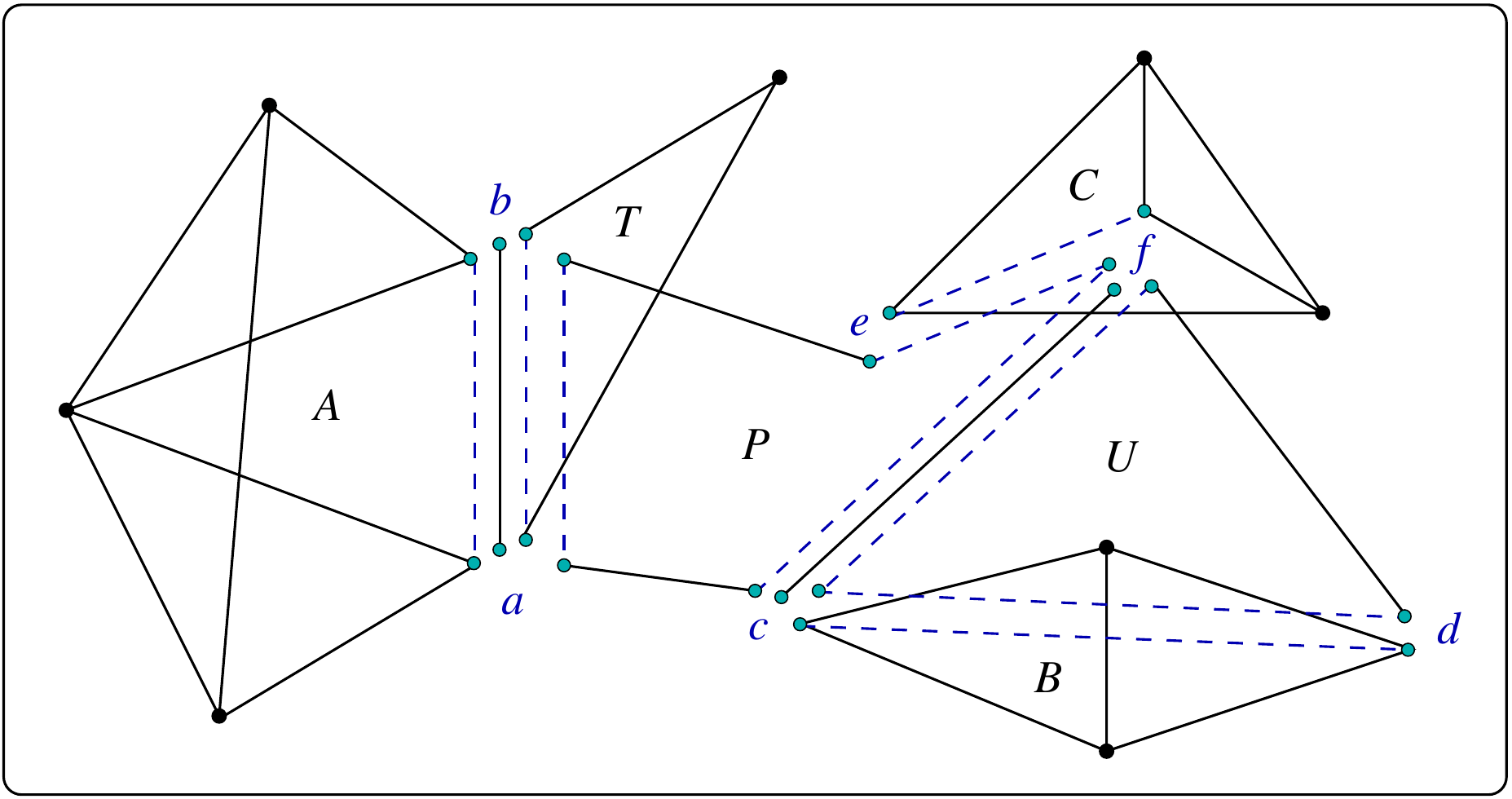}
\end{center}
\vspace{-4mm}
\caption{\small The 3-decomposition of $g$
\label{fig:tridecomp}}
\end{figure}

We also note whether the edge $xy$ is present or not in the original graph $g$. 
This has been done in Figure \ref{fig:tridecomp} for the separating pairs $\{a,b\}$, $\{c,d\}$, $\{c,f\} $ and $\{e,f\}$ and we see the 3-components appearing: the 3-connected components $A$, $B$, and $C$, and the polygons $P$, $T$ and $U$. The above dissection could also be performed for the separating pairs $\{a,e\}$ or $\{c,e\}$, for example, but that would cut the polygon $P$ into smaller polygons; so it is not done since the maximality of the polygonal components ensures the unicity of the decomposition. We refer the reader to the bibliography for more details. The essential separating pairs of $g$, ($\{a,b\}$, $\{c,d\}$, $\{c,f\} $ and $\{e,f\}$ in the example) will be referred to as the \emph{bonds} of the 3-decomposition. Hence a bond $\{x,y\}$ links together 2 or more 3-components, together with possibly the edge $xy$, with the exception of two polygons alone which is forbidden.

\begin{figure}[h] 
\begin{center}
\includegraphics[width=3.5in]{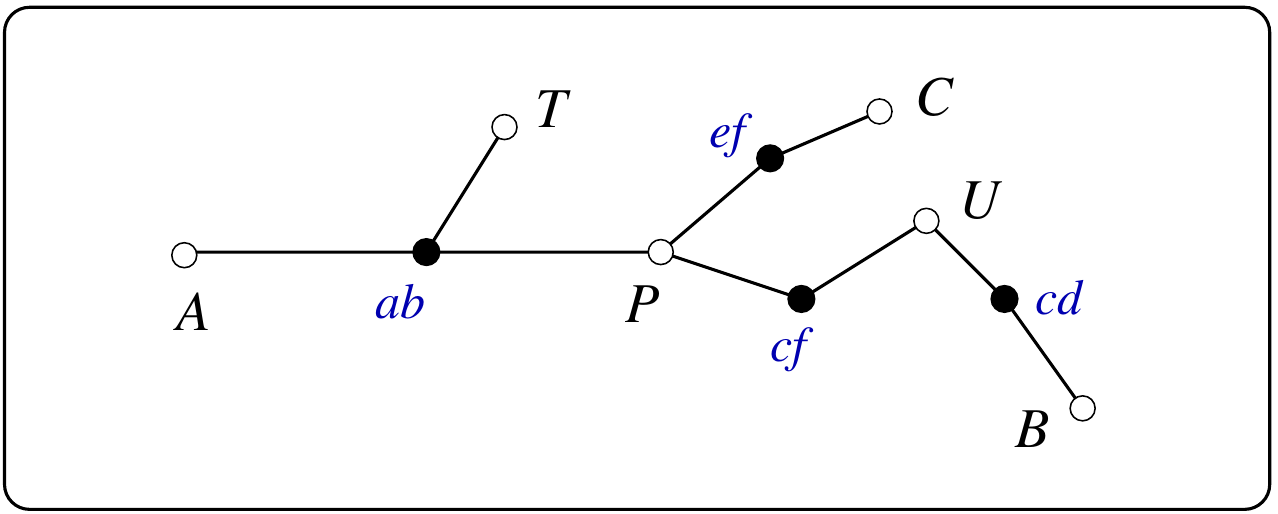}
\end{center}
\vspace{-4mm}
\caption{\small The $\tc$-tree of $g$
\label{fig:tcarbre}}
\end{figure}

 By analogy with the block-cutpoint tree $\bc(c)$ of a connected graph $c$, we deduce from this decomposition a bicolored tree $\tc(g)$ associated with any 2-connected graph $g$, whose white vertices are the 3-components of $g$ (3-connected components or polygons)  and whose black vertices are the separating pairs linking together these 3-components (the bonds). See Figure \ref{fig:tcarbre}.

Now let $\F$ be a given species of 3-connected graphs and $\B = \B_{\F}$ be the class of 2-connected graphs all of whose 3-connected components are in $\F$. Note: by convention,  $K_2$ is in $\B$.
%
We introduce the following classes of \emph{rooted} 
graphs in $\B$, relative to the concept of $\tc$-tree:

$\B^{\,\circ}$ denotes the class of all graphs $g \in \B$, rooted at a distinguished 3-component (3-connected component or polygon);

$\B^{\,\bullet}$ denotes the class of all graphs $g \in \B$ together with a distinguished bond;

$\B^{\,\circ-\bullet}$ denotes the class of all graphs $g \in \B$ with a distinguished pair of adjacent 3-component and bond. 

%
\begin{theorem} \textbf{\emph{(Dissymmetry Theorem for 2-connected graphs).}} \label{theo:dissym}
Let $\F$ be a given species of 3-connected graphs and $\B = \B_{\F}$ be the class of 2-connected graphs all of whose 3-connected components are in $\F$. We then have the following identity (species isomorphism):
\begin{equation} \label{eq:dissym}
\B^{\,\circ}\ +\ \B^{\,\bullet}\ = \ \B\, -\, K_2 \  +\ \B^{\,\circ-\bullet}\, .
\end{equation}
\end{theorem}
\textbf{Proof.} The proof uses the concept of \emph{center} of a tree. Notice that all the leaves of a $\tc$-tree are of the same color (white). This implies that its center is a vertex, black or white. 
Now a structure $s$ belonging to the left-hand side of 
(\ref{eq:dissym}) is a graph $g \in \B$ which is rooted at either a tricomponent or a bond of $g$, that is at a white or black vertex of $\tc(g)$.  It can happen that the rooting is performed right at the center. This is canonically equivalent
to doing nothing and is represented by the term $\B - K_2$ in the right-hand side of
(\ref{eq:dissym}). On the other hand, if the rooting is done at an off-center vertex, black or white, then there is a unique adjacent vertex of the other color in $\tc(g)$ which is located 
closer to the center, thus defining a unique $\B^{\,\circ-\bullet}$-structure. It is easily checked that this correspondence is bijective and independent of any labelling, giving the desired species isomorphism.
\hfill\rule{2mm}{2mm}

\medskip
Our next goal is to find closed form expressions for the species $\B^{\,\circ}$, $\B^{\,\bullet}$ and $\B^{\,\circ-\bullet}$.  This will be achieved using the operation of substitution of 01-networks for the edges of \emph{core} graphs, as explained in the next section.
%
%
\section{Operations on networks and their exponential generating functions.}
We first describe the exponential generating functions which are used for the labelled enumeration of graphs and networks and of related species, according to the number of edges.
For a species $\G$ of graphs 
the exponential generating function $\G(x,y)$, where the variable $y$ acts as an edge counter, is defined by 
\begin{eqnarray}  
\G(x,y)
=\sum_{n\ge 0}g_n(y)\frac{x^n}{n!}=\sum_{n\ge 0}\sum_{m\ge 0}g_{n,m}y^m\frac{x^n}{n!},
\end{eqnarray}
%
where $g_{n,m}$ is the number of graphs in $\G$ with $m$ edges over a given $n$-element set of vertices.  Similar definitions can be given for associated (for example rooted) species of graphs. For a species $\N$ of 01-networks, the exponential generating function $\N(x,y)$ is defined by
\begin{eqnarray}  
\N(x,y)=\sum_{n\ge 0}\nu_n(y)\frac{x^n}{n!}= \sum_{n\ge 0}\sum_{m\ge 0}\nu_{n,m}y^m\frac{x^n}{n!} ,
\end{eqnarray}
where $\nu_{n,m}$ is the number of 01-networks in $\N$ with $m$ edges over an $n$-element set of internal vertices.

We define an operator $\tau$ acting on 01-networks, $N \longmapsto \tau\cdot N$, that interchanges the poles $0$ and $1$. 
A class $\N$ of networks is called \textit{symmetric} if $N\in\N \Longrightarrow \tau\cdot N\in\N$.

%
If $\M$ is a species of networks not containing the edge $01$, then we denote by $y\,\M$ the class obtained by adding this edge to all the networks in $\M$. Observe that there are two distinct networks on the empty set, namely, the \emph{trivial network} $\nn$ consisting of two isolated vertices 0 and 1, and the one-edge network $y\nn$. 

Let $B$ be a given species of $2$-connected graphs, 
for example $B=\B_{all}$, the class of all $2$-connected graphs, $B=K_2$,
or $B=\B_P$, the class of all $2$-connected planar graphs.
We denote by $B^{(y)}$ the species of graphs obtained by selecting and removing an edge
from graphs in $B$ in all possible ways. 
Note that the removed edge is remembered so that
\begin{equation}
B^{(y)}(x,y)=\frac{\partial}{\partial y}B(x,y).
\end{equation}
If, moreover, the endpoints of the selected and removed edge are unlabelled and numbered as $0$ and $1$, in all possible ways, the resulting class of networks is denoted by $B_{0,1}$. For example, $(\B_{all})_{0,1}$ is the class of networks having non-adjacent poles, $(K_2)_{0,1} = \nn$, and the class of strongly planar networks can be expressed as $\N_P = (1+y)(\B_P)_{0,1} - \nn$, where the multiplication $y\cdot B_{0,1}$ corresponds to adding the edge $01$ to all networks in $B_{0,1}$. As another example,  $(K_5)_{0,1}$, is illustrated in Figure \ref{fig:couronne}. 
Relabelling the two poles yields the identity 
\begin{equation}
x^2B_{0,1}(x,y)=2\,B^{(y)}(x,y)
\end{equation}
so that
\begin{equation}
B_{0,1}(x,y)=\frac{2}{x^2}\frac{\partial}{\partial y}B(x,y).
\end{equation}
%

\subsection{Series composition}
\begin{defn}
\emph{
a) Let $M$ and $N$ be two non-trivial disjoint networks. 
The \emph{series composition} of $M$ followed by $N$, denoted by $M\cdot_sN$, is a 
vertex-rooted network whose underlying set is the union of the underlying sets of $M$ and $N$ plus an extra element. It is obtained by taking the graph union of $M$ and $N$, where the $1$-pole of $M$ is identified with the $0$-pole of $N$, and this \emph{connecting vertex} is labelled by the extra element and is the root of $M\cdot_sN$. See Figure \ref{fig:seriescomp}.\\
b) The underlying (unrooted) network of a series composition is called an $s$-\emph{network}.\\
c) The \emph{series composition} $\M\cdot_s\N$ of two species of non-trivial networks $\M$ and $\N$  
is the class obtained by taking all series compositions $M\cdot_s N$ of networks with $M\in\M$ and $N\in\N$. 
} 
\end{defn}
\begin{figure}[h] \label{fig:seriescomp}
\begin{center} \includegraphics[width=5.0in]{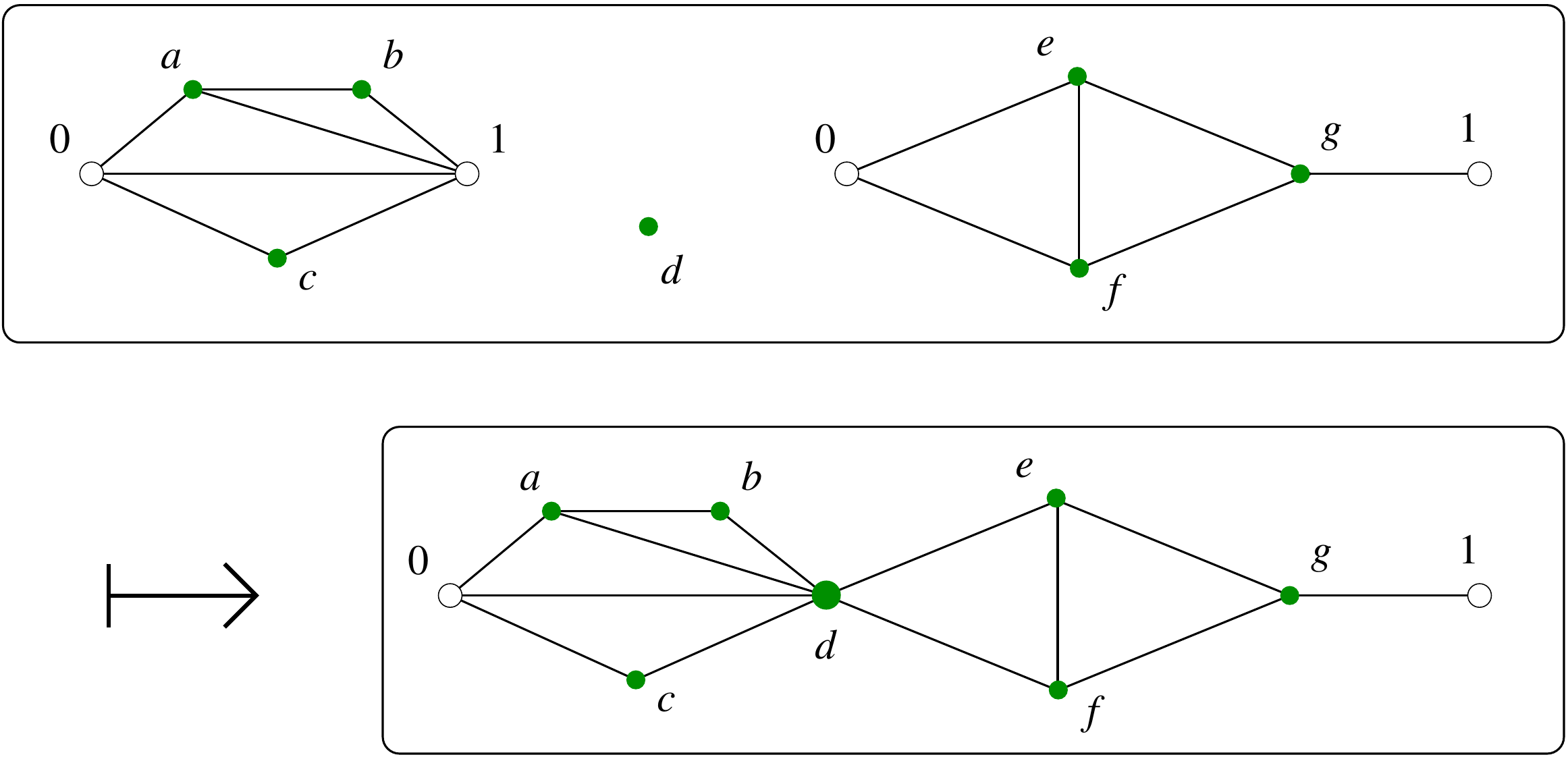}
\end{center}
\vspace{-2mm}
\caption{Series composition of networks\label{fig:seriescomp}}
\end{figure}
The species $\M\cdot_s\N$ can be expressed as the species product $\M\, X\,\N$, where the factor $X$ corresponds to the connecting vertex, and we have
\begin{equation} \label{eq:series}
(\M\cdot_s\N)(x,y)=x\M(x,y)\N(x,y).
\end{equation}
If for any $\M\cdot_s\N$-structure the two components $M\in\M$ and $N\in\N$ are uniquely determined by the resulting network (for example if no $M\in\M$ is an $s$-network), then the series composition is called \emph{canonical} and the rooting of the connecting vertex can be neglected.  We say that a species of networks $\D$ is \emph{closed under series composition and decomposition} if for any $s$-network $R$, $R$ is in $\D$ if and only if each individual factor of $R$ is in $\D$.
\begin{propos} \label{prop:s-networks}
Let $\D$ be a species of networks which is closed under series composition and decomposition and let $\S$ denote the class of $s$-networks in $\D$. Then we have
\begin{equation} \label{eq:SDsD}
\S= (\D-\S)\cdot_s\D
\end{equation}
the series composition being canonical, and also
\begin{equation} \label{eq:SD2XD}
\S = 
\frac{X\D^2}{1+X\D}.
\end{equation}
\end{propos}
\proof  Any $s$-network in $\D$ can be decomposed uniquely into a first network which is not an $s$-network followed by an arbitrary network of $\D$, whence (\ref{eq:SDsD}). We also have $\S=(\D-\S)\,X\,\D$ and solving for $\S$ yields (\ref{eq:SD2XD}).   
\hfill\rule{2mm}{2mm}
\begin{corol}
Under the hypothesis of Proposition \ref{prop:s-networks}, we have, for the exponential generating function,
\begin{equation} \label{eq:Sxy}
\S(x,y) = \frac{x\D^2(x,y)}{1+x\D(x,y)}.
\end{equation}
\end{corol}
%
\textbf{Remark.}  
Iterating the idea behind the decomposition (\ref{eq:SDsD}), one has the more symmetric canonical decomposition
\begin{equation} \label{eq:SDmoinsS}
\S =\, (\D-\S)\cdot_s(\D-\S)\; +\; (\D-\S)\cdot_s(\D-\S)\cdot_s(\D-\S)\; + 
\cdots. 
\end{equation}
%

\subsection{Parallel composition}
\begin{figure}[h]
\begin{center} \includegraphics[width=4.5in]{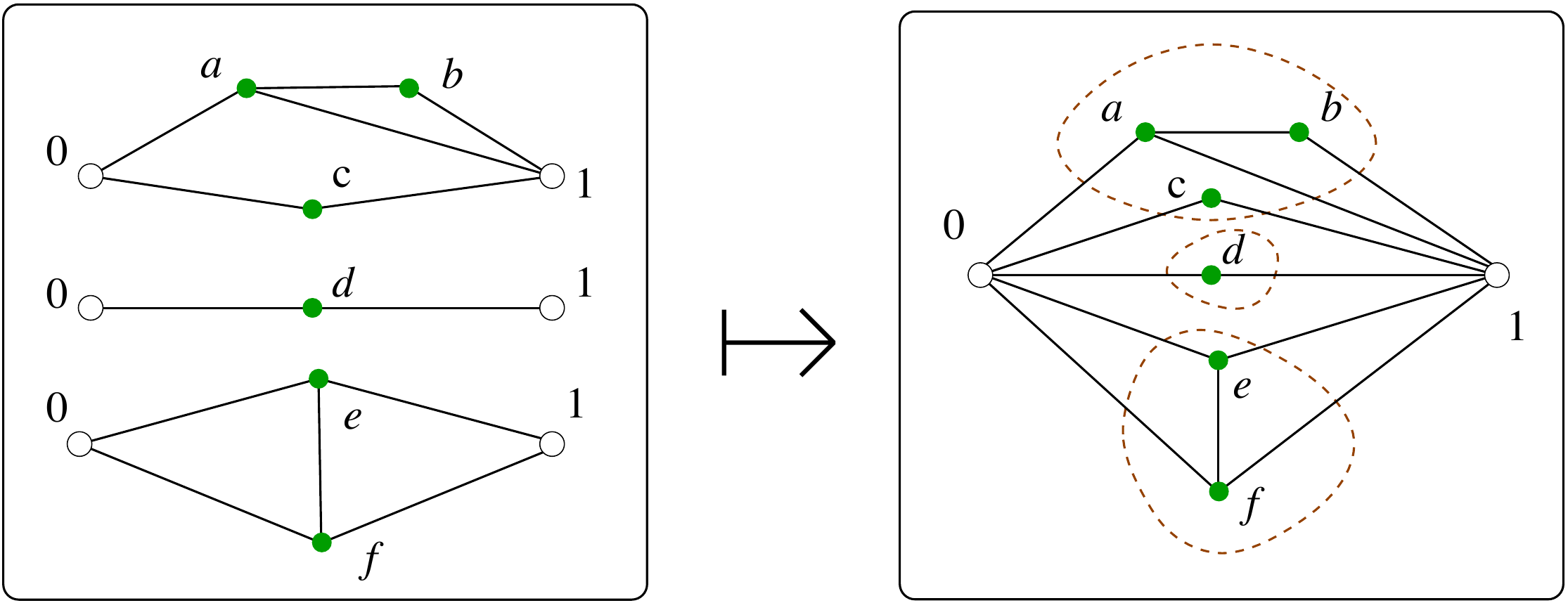}
\end{center}
\vspace{-4mm}
\caption{Parallel composition of networks \label{fig:parallelcomp}}
\end{figure}

\begin{defn} \label{def:parallel}
\emph{
a) Let $N$ be a finite set of disjoint non-trivial 01-networks having non-adjacent poles. The \emph{parallel composition} of $N$ is the \emph{partitioned} 01-network obtained by taking the union of the graphs in $N$, where all $0$-poles are merged into one $0$-pole, and similarly for the $1$-poles; the partition 
of the internal vertices into those of the networks of $N$ is part of the structure. 
An example is given in Figure \ref{fig:parallelcomp}.
By convention, the parallel composition of an empty set of networks is the trivial network $\nn$, while the parallel composition of one network is the network itself. \\
b) The underlying (unpartitioned) network of a parallel composition of at least two networks is called a $p$-\emph{network}. Any network having adjacent poles is also considered as a $p$-network. \\
c) If $\N$ is a species of non-trivial 01-networks having non-adjacent poles, then the \emph{parallel composition} of $\N$ is defined as the class of all parallel compositions of finite sets of disjoint networks in $\N$.  This operation is denoted by $E(\N)$, the ordinary composition of species, where $E$ denotes the species of sets, since a parallel composition can be seen as an assembly of networks, with all 0-poles identified and also all 1-poles. As mentioned above, $\nn=E_0(\N)$ and $\N = E_1(\N)$.
}
\end{defn}

If each network in a class $\M$ can be viewed unambiguously as a 
parallel composition of networks in $\N$, then we say that the parallel composition is called \emph{canonical}. This happens if and only if none of the networks in $\N$ is a $p$-network. Then we can write $\M=E(\N)$ and for the exponential generating functions, we have 
\begin{equation} \label{eq:parallelxy}
\M(x,y)=\exp(\N(x,y)).
\end{equation}
%
%
\subsection{The $\uparrow$-composition (substitution of networks for edges)}
\begin{defn}
\emph{
Let $\M$ be a species of graphs (or networks) and $\N$ be a symmetric species of non-trivial networks.  We denote by $\T=\M \uparrow \N$ the class of pairs 
$(M,T)$ such that
\begin{enumerate}
\item
the graph (or network) $M$ (called the \emph{core}) is in $\M$, 
\item
there exists a family $\{N_e\}$ of networks in $\N$ (called the \emph{components}) such that the graph $T$ can be obtained from $M$ by substituting  $N_e$ for each edge $e$ of $M$, the poles of $N_e$ being identified with the extremities of $e$. 
\end{enumerate}}
\end{defn}
\begin{figure}[h] 
\begin{center} \includegraphics[height=2.5in]{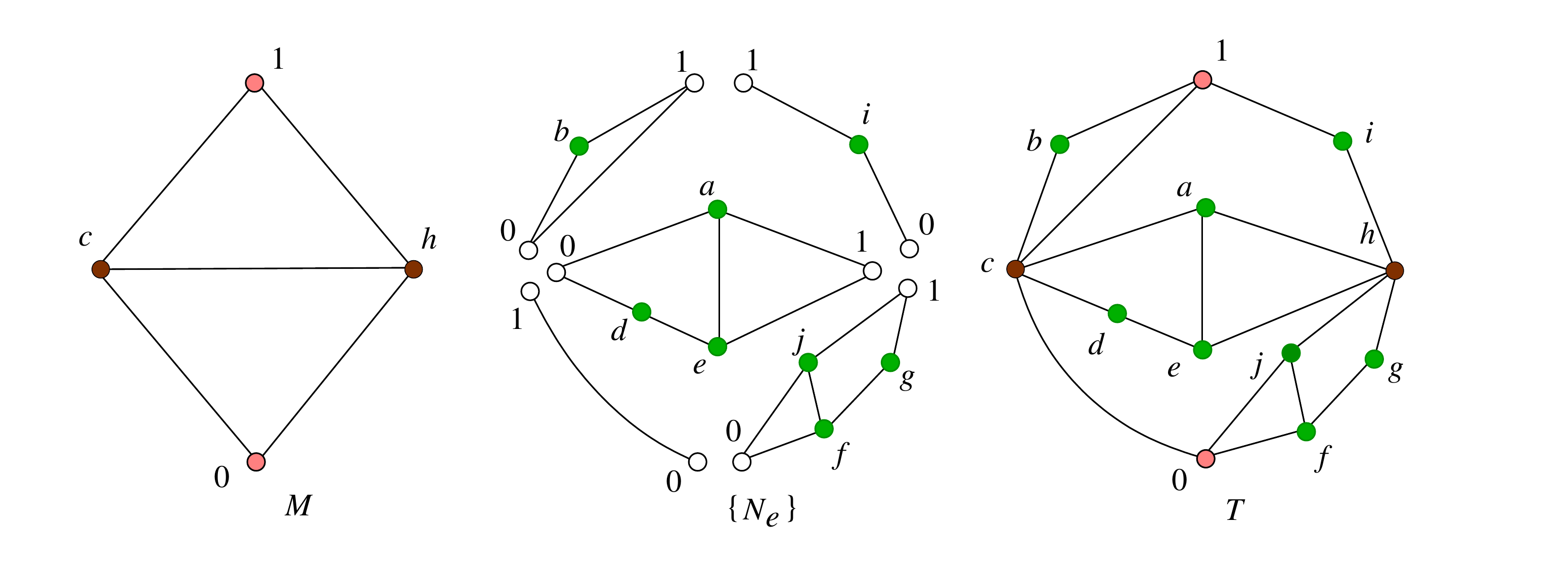}
\end{center}
\vspace{-8mm}
\caption{Example of a $(\M\uparrow \N)$-structure $(M,T)$, with $\M = K_4\backslash e$ \label{fig:exemple3}}
\end{figure}

An example of an $(\M \uparrow \N)$-structure $(M,T)$ is presented in Figure~\ref{fig:exemple3}, where $\M=K_4\backslash e$ and $\N$ is the class of all networks. Notice that if the network $N_e$ which is substituted for the edge $e=xy$ is not the one-edge network $y\nn$, then the pair of vertices $\{x,y\}$ is a separating pair of $T$.

As another example,  take $\G=K_2$ and let $\N$ be any symmetric species of networks. Then $K_2\uparrow \N$ consists of rooted graphs obtained from networks in $\N$ by labelling the two poles 0 and 1, the rooting being at this pair of vertices.
\begin{propos}
Let $\M$ be a species of graphs (or networks) and $\N$ be any symmetric species of non-trivial networks. Then we have, for the labelled enumeration,
\begin{equation} \label{eq:MdexNxy}
(\M \uparrow \N)(x,y) = \M(x,\N(x,y)).
\end{equation}
\end{propos}
\proof See \cite{Timothy} or \cite{GLL}. \hfill\rule{2mm}{2mm}

\medskip
The composition $\M \uparrow \N$ is called {\it canonical} if for any structure $(M,T)\in \M \uparrow \N$ the core  $M\in \M$ is uniquely determined by the graph (or network) $T$. In this case, we can identify $\M \uparrow \N$ with the species of resulting graphs (or networks) $T$.
%

%
\subsection{Functional equations}

Let $\F$ be a species of $3$-connected graphs and 
let $\B=\B_\F$ (resp. $\D=\D_\F$) denote the class of all 2-connected graphs (resp. non-trivial networks) all of whose 3-connected components are in $\F$, where the 3-components of a network $N$ are defined by applying the 3-decomposition to the graph $N\cup01$. Then
\begin{equation} \label{eq:DB01}
\D = (1+y)\B_{0,1} - \nn,
\end{equation}
where multiplication by $y$ corresponds to adding the edge $01$ to all networks in $\B_{0,1}$.

One example of a non-canonical $\uparrow$-composition is given by $\F\uparrow\D$ which represents the species of 2-connected graphs in $\B$ rooted at some 3-connected component. By contrast, as stated in Theorem \ref{theo:reseaux} below, any composition of the form $\F_{0,1}\uparrow\N$ is canonical. 

Let $\S$ denote the subclass of $\D$ consisting of $s$-networks. By virtue of Proposition \ref{prop:s-networks}, we have
\begin{equation} \label{eq:SXD2D}
\S = \frac{X\D^2}{1+X\D}.
\end{equation}

Parts a) and b) of the next theorem can be seen as a species form of Trakhtenbrot's decomposition theorem \cite{Trakh} which was originally stated and proved for networks in which parallel edges are allowed.  See \cite{Timothy, Timothyunlabelled} for a precise statement of Trakhtenbrot's Theorem. A proof in English is available from the fourth author on request.
%
%
%
\begin{thm} \label{theo:reseaux}
Let $\F$ be a species of $3$-connected graphs and $\D=\D_\F$ be the corresponding class of $01$-networks associated to $\F$.  Then:\\
\emph{a)} For any (symmetric) species $\N$ of non-trivial networks, the composition $\F_{0,1}\uparrow\N$ is canonical.\\
\emph{b)} Let $\H$ denote the subclass of $h$-networks, defined by  
\begin{equation} \label{eq:HF01D}
\H = \F_{0,1}\uparrow\D.
\end{equation}
Then we have 
\begin{equation} \label{eq:DSPH}
\D=\S+\P+\H,
\end{equation}
where $\P$ denotes the subclass of $p$-networks of $\D$ \emph{(see Definition \ref{def:parallel} b))}, which satisfies 
\begin{equation} \label{eq:PdeHS}
\P = (1+y)E_{\geq2}(\H+\S)+ y(\H+\S) + y\nn.
\end{equation}
\emph{c)} The species $\D$ is characterized by the functional relation
\begin{equation} \label{eq:RdeR}
\D = (1+y)E(\F_{0,1}\uparrow\D + \frac{X\D^2}{1+X\D})- \mathrm{\nn}. 
\end{equation}
\end{thm}
\proof 
As mentioned above, parts a) and b) are essentially a reformulation of Trakhtenbrot's decomposition theorem for networks, where parallel edges are not allowed.  Intuitively, for a), observe that in any $\F_{0,1}$-network $M$, the poles are non-adjacent. Hence, the same will be true for any network $T$ arising from a $\F_{0,1}\uparrow\N$-structure $(M,T)$. Applying the 3-decomposition of 2-connected graphs of Section 2 to the graph $T\cup01$, we can see that there is a unique 3-connected component in $T$ containing the vertices 0 and 1, namely $M$ itself.\\
For b), it is easy to see whether a network $N$ having non-adjacent poles is an $s$-network (the graph is not 2-connected) or a $p$-network (the poles form a separating pair) and, otherwise, that $N$ is in fact of the form $\F_{0,1}\uparrow\D$, again using the 3-decomposition of 2-connected graphs. If the network has adjacent poles, then it must be of the form $y(\nn + \H+\S + E_{\geq2}(\H+\S)= yE(\H+\S))$ so that (\ref{eq:DSPH}) and (\ref{eq:PdeHS}) are satisfied.\\
c) Putting (\ref{eq:DSPH}) and (\ref{eq:PdeHS}) together, we find that
\begin{eqnarray} 
\D&=&\H+\S+ y\nn + y(\H+\S) + (1+y)E_{\geq2}(\H+\S) \nonumber  \\
&=& (1+y)E(\H + \S) - \mathrm{\nn}  \label{eq:DEHS}
\end{eqnarray}
and (\ref{eq:RdeR}) follows from (\ref{eq:DEHS}), (\ref{eq:HF01D}) and (\ref{eq:SXD2D}).  
\hfill\rule{2mm}{2mm}

\medskip
We are now in position to express the three rootings $\B^{\,\circ}$, $\B^{\,\bullet}$ and $ \B^{\,\circ-\bullet}$ of the species $\B=\B_\F$ associated with a given species of 3-connected graphs $\F$, which occur in the Dissymmetry Theorem for 2-connected graphs (Theorem \ref{theo:dissym}), in terms of the corresponding class $\D=\D_\F$ of 01-networks. Let $\C$ denote the species of polygons, that is of (unoriented) cycles of length $\geq3$.

\begin{theorem} \label{theo:rootingsD}
We have the following identities:
\begin{eqnarray}
\B^{\,\circ} & = & \F \uparrow \D \ +\  \C\uparrow (\D-\S),\label{eq:Brond} \\
\B^{\,\bullet}   %
& = & K_2 \uparrow \left((1+y)E^{\phantom{2}}_{\geq2}(\H+\S) - E_{2}(\S) \right)\label{eq:Bpoint}\\
& = & K_2 \uparrow \left(\D - (1+y)(\H + \S) - y\nn - E^{\phantom{2}}_2(\S)\right), \label{eq:Bpointbis}\\
\B^{\,\circ-\bullet} 
& = & K_2 \uparrow \left((\H + \S)(\D - y\nn) -\S^2 \right), \label{eq:Brondpoint}
\end{eqnarray}
where $\D$ is characterized by equation $(\ref{eq:RdeR})$,  $\S= \frac{X\D^2}{1+X\D}$ and $\H=\F_{0,1}\uparrow \D$.
\end{theorem}
\proof Recall that
$\B^{\,\circ}$ denotes the class of all graphs $g \in \B$ rooted at a distinguished 3-component. If the distinguished component $C$ is a 3-connected graph, then $g$ is  obtained by replacing each edge of $C$ by a network in $\D$. Otherwise, $C$ is a polygon and $g$ will be obtained by replacing each edge of $C$ by a network in $\D-\S$, by virtue of the maximality of polygonal 3-components. This establishes (\ref{eq:Brond}).

Also recall that the $K_2\uparrow$ operator transforms networks into rooted graphs by labelling the two poles 0 and 1, the rooting being at this pair of vertices. 
Now $\B^{\,\bullet}$ denotes the class of all graphs $g \in \B$ together with a distinguished bond. This separating pair $\{a,b\}$ decomposes g into two or more pieces which can be seen as either $h$-networks or $s$-networks. This yields the term $(1+y)E^{\phantom{\sum}}_{\geq2}(\H+\S)$, the factor $(1+y)$ accounting for the possibility that the vertices $a$ and $b$ be adjacent. However, the case of a non-adjacent separating pair joining two $s$-networks has to be excluded since this would imply decomposing a polygon into two smaller ones, which is prohibited. Hence (\ref{eq:Bpoint}). Formula (\ref{eq:Bpointbis}) then follows easily, using (\ref{eq:DEHS}).

The proof of (\ref{eq:Brondpoint}) is similar, the difference being that now one of the separated components is also distinguished. Details are left to the reader.
%
\hfill\rule{2mm}{2mm}
\begin{corol} [Explicit form of the Dissymmetry Theorem] 
\label{theo:dissymBis}
Let $\F$ be a given spe\-cies of 3-connected graphs and $\B = \B_{\F}$ be the species of 2-connected graphs all of whose 3-connected components are in $\F$. Also let  $\D=\D_\F$ denote the corresponding species of 01-networks. We then have the following identity: 
\begin{eqnarray} 
\B & = & \F \uparrow \D \ +\  \C\uparrow (\D-\S) \nonumber  \\ 
& &  \phantom{M}+ K_2 \uparrow \left(\D -(\H+\S)(\D+1) - E_2(\S) + \S^2 \right), \label{eq:dissymBis}
\end{eqnarray}
where $\D$ is characterized by equation $(\ref{eq:RdeR})$,  $\S= \frac{X\D^2}{1+X\D}$ and $\H=\F_{0,1}\uparrow \D$.
\end{corol}
%
\section{Series techniques for the unlabelled enumeration of graphs and networks}
Traditionally, two generating series are used for the unlabelled enumeration of structures: the ordinary (tilde) generating function and the cycle index series. These are now reviewed in the context of graphs and networks where the number of vertices and the number of edges are taken into account and where a variant of the cycle index series is necessary when dealing with the $\uparrow$-composition, the substitution of networks for edges. This variant, the \emph{edge index series}, is introduced in \cite{Timothyunlabelled} and called Walsh index series in \cite{GLL3}. Detailed proofs of most of their main properties can be found in \cite{GLL3}.

For a species $\M$ of (possibly rooted) graphs or networks,
the ordinary (\emph{tilde}) generating function $\widetilde{\M}(x,y)=\M\tilde{\ }(x,y)$ is defined as follows: 
\begin{eqnarray}
\widetilde{\M}(x,y)
=\sum_{n\ge 0}\tilde{\mu}_n(y)x^n=\sum_{n\ge 0}\sum_{m\ge 0}\tilde{\mu}_{n,m}y^mx^n,
\end{eqnarray}
where $\tilde{\mu}_{n,m}$ is the number of isomorphism classes of graphs (resp. networks) in $\M$ having $n$ vertices (resp. internal vertices) and $m$ edges.


Let $\G$ be a species of graphs and let $\N$ be a species of networks. The three edge index series
$W_\G({\bf a}; {\bf b}; {\bf c})$, $W^+_\N({\bf a}; {\bf b}; {\bf c})$ and $W^-_\N({\bf a}; {\bf b}; {\bf c})$, in variables ${\bf a}=(a_1,a_2,\ldots)$, ${\bf b}=(b_1,b_2,\ldots)$ and ${\bf c}=(c_1,c_2,\ldots)$ are defined in what follows.

\begin{figure}[h]
\begin{center}
	\includegraphics[width=4.0in]{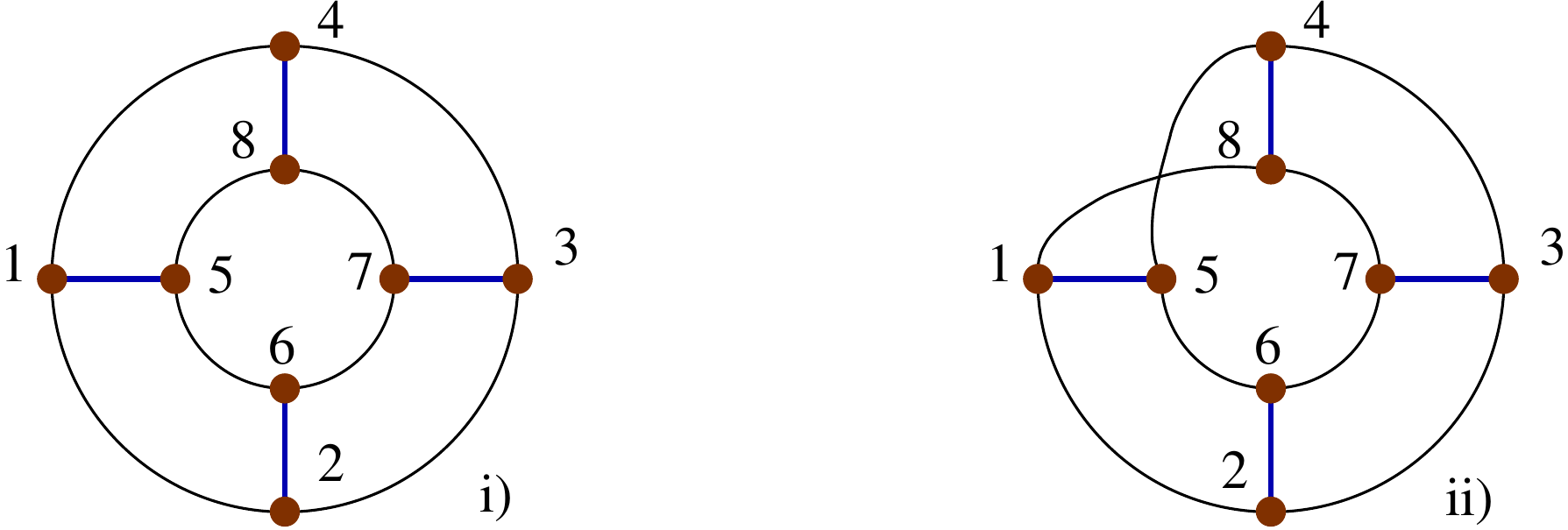}
\end{center}
\vspace{-5mm}
	\caption{(i) Cylindrical edge cycle, (ii) M\"obius edge cycle.}
	\label{fig:cylmoe}
\end{figure}

Let $G=(V(G),E(G))$ be a graph in $\G$. A permutation $\sigma$ of $V(G)$ that is an automorphism of the graph $G$ induces a permutation $\sigma^{(2)}$ of the set $E(G)$ of edges whose cycles are of two possible sorts: if $c$ is a cycle of $\sigma^{(2)}$ of length $l$, then either $\sigma^l(a)=a$ and $\sigma^l(b)=b$ for each edge $e=ab$ of $c$ (a \textit{cylindrical} edge cycle), or else $\sigma^l(a)=b$ and $\sigma^l(b)=a$ for each edge $e=ab$ of $c$ (a \textit{M\"obius} edge cycle). For example, the automorphism $\sigma=(1,2,3,4)(5,6,7,8)$ of the graph of Figure~\ref{fig:cylmoe} (i) induces the cylindrical edge cycle 
$(15,26,37,48)$, and the automorphism $\sigma=(1,2,3,4,5,6,7,8)$ of the graph of Figure~\ref{fig:cylmoe} (ii) induces the M\"obius edge cycle
$(15,26,37,48)$.

For an automorphism $\sigma\in\Aut(G)$ of $G$, denote by $\sigma_k$ the number of cycles of length $k$ of $\sigma$, by $\mathrm{cyl}_k(G,\sigma)$ the number of cylindrical edge cycles of length $k$,
and by $\mathrm{m\ddot{o}b}_k(G,\sigma)$ the number of M\"obius edge cycles of length $k$ induced by $\sigma$ in $G$.
Given a graph $G\in \G$ and an automorphism $\sigma$ of $G$, the \textit{weight} $w(G,\sigma)$ of such a structure is the following cycle index monomial:
\begin{eqnarray}
w(G,\sigma)=a_1^{\sigma_1}a_2^{\sigma_2}\cdots
b_1^{\mathrm{cyl}_1(G,\sigma)}b_2^{\mathrm{cyl}_2(G,\sigma)}\cdots c_1^{\mathrm{m\ddot{o}b}_1(G,\sigma)}c_2^{\mathrm{m\ddot{o}b}_2(G,\sigma)}\cdots\,.
\label{form:weightGsigma}
\end{eqnarray}
The \emph{edge index series} $W_\G({\bf a}; {\bf b}; {\bf c})$ of $\G$ is defined as 
\begin{eqnarray}
W_\G({\bf a}; {\bf b}; {\bf c}) = \sum_{G\underline{\in} \mathrm{Typ}(\G)} \frac{1}{|\Aut(G)|} \sum_{\sigma\in \Aut(G)} w(G,\sigma),
\label{form:Wabc}
\end{eqnarray}
where the notation $G\underline{\in} \mathrm{Typ}(\G)$ means that the summation should be taken over a set of representatives $G$ of the isomorphism classes of graphs in $\G$.

\medskip  
\noindent
\textbf{Examples.}  1.  The edge index series of the species $K_2$ is given by 
\begin{equation}
W_{K_2}  =  \frac{1}{2}(a_1^2b_1 + a_2c_1).
\end{equation}
2. The edge index series of the species $\C_n$ of (unoriented) cycles of length $n$ is a refinement of the usual cycle index polynomial for $\C_n$. It is given by
\begin{equation}
W_{C_n}({\bf a}; {\bf b}; {\bf c})=\frac{1}{2n}\sum_{d|n}\phi(d)a_d^{\frac{n}{d}}b_d^{\frac{n}{d}} + 
\frac{1}{2}\left\{ \begin{array}{ll} 
a_1a_2^{\frac{n-1}{2}}b_2^{\frac{n-1}{2}}c_1, & \mbox{\emph \,$n$ }\mathrm{odd}\\ 
\frac{1}{2}(a_2^{\frac{n}{2}}b_2^{\frac{n-2}{2}}c_1^2 + a_1^2a_2^\frac{n-2}{2}b_2^{\frac{n}{2}}), & \mbox{\emph \,$n$ }\mathrm{even} 
\end{array}
\right.                                          
\label{eq:WCn}
\end{equation}
where $\phi$ is the Euler $\phi$-function.  See  \cite{Timothyunlabelled} and \cite{GLL3} which contains a typo in formula (43). 
By summing over $n\geq3$, we obtain the edge index series of $\C$. 
The result is
\begin{eqnarray}
W_{\C}  &=&  \frac{1}{2}\sum_{d\geq1}\frac{\phi(d)}{d}\log\frac{1}{1-a_d b_d}-\frac{1}{2}a_1b_1-\frac{1}{4}a_1^2b_1^2-\frac{1}{4}a_2b_2  \nonumber\\
&&\ +\frac{1}{4}(2a_1c_1 + a_2c_1^2 + a_1^2b_2)\frac{a_2b_2}{1-a_2b_2}.\label{eq:WC}
\end{eqnarray}

\medskip
Note that any isomorphism of networks $\varphi : N \tilde{\longrightarrow} N^\prime$ is assumed to be pole-preserving, i.e. $\varphi(0)=0$ and $\varphi(1)=1$. In particular, any automorphism of a network $N$ should be pole-preserving. 
It will be necessary to consider the subclass $\N_\tau$ of $\N$ consisting of  $\tau$-\textit{symmetric networks}, i.e. 
\begin{equation}
\N_\tau=\{ N\in \N\ |\ \tau\cdot N\simeq N \}.
\label{form:A+}
\end{equation}

Let $U$ be the underlying set of a network $N$ and suppose that $\sigma$ is in $S[U]$, i.e. $\sigma$ is a permutation of $U$. We can extend $\sigma$ to permutations on $U\cup\{0,1\}$, $\sigma^+=(0)(1)\sigma$
and $\sigma^{-}=(0,1)\sigma$; in other words,  $\sigma^+$ preserves the poles and $\sigma^{-}$ exchanges them. 
For any network $N\in\N$, denote by $\hat{N}$ the corresponding graph on $U\cup \{0,1\}$. Then we introduce the notation
\begin{equation}
\Aut^+(N)=\{ \sigma\in S[U]\ |\ \sigma^+\in \Aut(\hat{N}) \}
\label{form:A+}
\end{equation}
and 
\begin{equation}
\Aut^-(N)=\{ \sigma\in S[U]\ |\ \sigma^-\in \Aut(\hat{N}) \}.
\label{form:A-}
\end{equation}
In other words, a $\sigma^+$ in (\ref{form:A+}) is a pole-preserving graph automorphism, i.e. a network automorphism, while a $\sigma^-$ in (\ref{form:A-}) is a pole-reversing \emph{graph} automorphism.
Notice that $\Aut^+(N)=\Aut(N)$ and that if $\Aut^-(N)$ is not empty, then 
$|\Aut^-(N)|=|\Aut^+(N)|$. This can be seen by using the composition of automorphisms.
For $N\in\N$ and $\sigma\in \Aut^+(N)$, we assign the weight
\begin{equation}
w(N,\sigma)=\frac{w(\hat{N},\sigma^+)}{a_1^2},
\label{form:weightN+} 
\end{equation}
where the second $w$ is defined by (\ref{form:weightGsigma}), and for $N\in \N$ and $\sigma\in \Aut^-(N)$, we set 
\begin{equation}
w(N,\sigma)=\frac{w(\hat{N},\sigma^-)}{a_2}.
\label{form:weightN-}
\end{equation}
In other words, only the internal vertex cycles are accounted for.
Then, for a species $\N$ of networks, the following two \textit{edge index series} are defined by
\begin{eqnarray} 
W^+_\N({\bf a}; {\bf b}; {\bf c})= \sum_{N\underline{\in} \mathrm{Typ}(\N)} \frac{1}{|\Aut^+(N)|} \sum_{\sigma\in \Aut^+(N)} w(N,\sigma),
\end{eqnarray}
\begin{eqnarray} 
W^-_\N({\bf a}; {\bf b}; {\bf c})= \sum_{N\underline{\in} \mathrm{Typ}(\N_\tau)} \frac{1}{|\Aut^-(N)|} \sum_{\sigma\in \Aut^-(N)} w(N,\sigma).
\end{eqnarray}

As the next proposition shows, the edge index series contain all the enumerative (labelled and unlabelled) information.
\begin{propos} [\cite{GLL3,Timothyunlabelled}] \label{propos:specDeW} Let $\G$ be a species of graphs and $\N$ be a species of networks. Then the following series identities hold:
\begin{eqnarray}
\G(x,y) & = & W_\G(x,0,0,\ldots; y,y^2,y^3,\ldots;  y,y^2,y^3,\ldots),\\
\widetilde{\G}(x,y) &=& W_\G(x,x^2,x^3,\ldots; y,y^2,y^3,\ldots; y,y^2,y^3,\ldots),\\
\N(x,y) & = & W_{\N}^+(x,0,0,\ldots; y,y^2,y^3,\ldots;  y,y^2,y^3,\ldots),\\
\widetilde{\N}(x,y) &=& W_{\N}^+(x,x^2,x^3,\ldots; y,y^2,y^3,\ldots; y,y^2,y^3,\ldots),\\
\N_{\tau}(x,y) & = & W_{\N}^-(x,0,0,\ldots;  y,y^2,y^3,\ldots;  y,y^2,y^3,\ldots),\\
\widetilde{\N}_{\tau}(x,y) &=& W_{\N}^-(x,x^2,x^3,\ldots; y,y^2,y^3,\ldots; y,y^2,y^3,\ldots).
\end{eqnarray}
\end{propos}
\hfill \rule{2mm}{2mm}

\medskip
Another description of the edge index series is very useful for understanding them and for establishing their properties.  It consists of expressions which involve exponential generating functions of labelled enumeration. These are recalled 
from Section 6 of \cite{GLL3}.

Following an idea of Joyal \cite{Joyal}, we introduce the auxiliary weighted species $\G^\aut=\G^\aut_w$. For any finite set $U$ (of vertices), $\G^\aut[U]$ is defined as the set of graphs in $\G[U]$ equipped with an automorphism $\sigma$, i.e. 
\[
\G^\aut[U]=\{ (G,\sigma)\ |\ G\in \G[U], \sigma\in S[U] : \sigma\cdot G=G \},
\]
where $S[U]$ is the set of all permutations of $U$. The relabelling rule of $\G^\aut$-structures along a bijection $\beta : U\ \tilde{\longrightarrow}\ U^\prime$ is defined as follows:
\[
\beta \cdot (G,\sigma) = (\beta \cdot G,\ \beta {\circ} \sigma {\circ} \beta^{-1}),
\]
where $\beta \cdot G$ is the graph obtained from $G$ by relabelling along $\beta$ and the composition $\circ$ is taken from right to left. It is easy to verify that $\G^\aut_w$ is a well-defined weighted species, where the weight function $w(G,\sigma)$ is the cycle index monomial defined by (\ref{form:weightGsigma}). Recall that $|\G^\aut[n]|_w$ denotes the total weight of $\G^\aut_w$-structures over the vertex set $[n]:=\{1,2,\ldots,n\}$, i.e. 
\begin{eqnarray*}
|\G^\aut[n]|_w=\sum_{(G,\sigma)\in \G^\aut_w[n]} w(G,\sigma).
\label{form:WabcAut}
\end{eqnarray*}
\begin{propos}[\cite{GLL3}] 
Using the exponential generating function of labelled $\G^\aut_w$-structures, we have
\begin{equation}
W_\G({\bf a}; {\bf b}; {\bf c}) = \sum_{n\ge 0}\frac{1}{n!}|\G^\aut[n]|_w 
= \G^\aut_w(x)\,|_{x=1}. 
\end{equation}
\end{propos}
\proof  The proof follows from the fact that the number of distinct graphs on $[n]$ obtained by relabelling a given graph $G$ with $n$ vertices is given by $\frac{n!}{|\mathrm{Aut}(G)|}$.  
\hfill \rule{2mm}{2mm}

\medskip
A similar approach can be used for the edge 
index series $W_\N^+$ and $W_\N^-$ of a given species of $2$-pole networks $\N$. We introduce the sets 
$$\N^+[U]=\{ (N,\sigma)\ |\ N\in \N[U],\ \sigma\in \mathrm{Aut}^+(N) \}$$ 
and 
$$\N^-[U]=\{ (N,\sigma)\ |\ N\in \N[U],\ \sigma\in \mathrm{Aut}^-(N) \},$$
where $\Aut^+(N)$ and $\Aut^-(N)$ are defined by (\ref{form:A+}) and (\ref{form:A-}), respectively.
Then, using the weight functions given by (\ref{form:weightN+}) and (\ref{form:weightN-}),
$\N^+_w$ and $\N^-_w$ are weighted species whose labelled enumerations yield by specialization the series $W^+_\N$ and  $W^-_\N$.

\begin{propos}[\cite{GLL3}] 
For a species of networks $\N$, the edge index series $W^+_\N$ and $W^-_\N$ can be expressed by the formulas
\begin{equation}
W^+_\N({\bf a}; {\bf b}; {\bf c}) = \N^+_w(x) |_{x=1}\,, \ \ \ \ 
 W^-_\N({\bf a}; {\bf b}; {\bf c}) =  \N^-_w(x) |_{x=1}\,. 
\label{form:35}
\end{equation}
%
\end{propos}
\hfill \rule{2mm}{2mm}

\smallskip
In order to describe the edge index series of an $\uparrow$-composition, we introduce the following  familiar plethystic notation. For any series of edge index type $f({\bf a}; {\bf b}; {\bf c})$ and any integer $k\geq1$, we set
%
\begin{equation} \label{eq:plethisme}
f_k=f_k({\bf a}; {\bf b}; {\bf c}) = f(a_k, a_{2k},a_{3k},\ldots; b_k, b_{2k},b_{3k},\ldots; c_k, c_{2k},c_{3k},\ldots). 
\end{equation}
Morevover, for any series $\ell=\ell({\bf a}; {\bf b}; {\bf c})$, $f=f({\bf a}; {\bf b}; {\bf c})$, $g=g({\bf a}; {\bf b}; {\bf c})$, $h=h({\bf a}; {\bf b}; {\bf c})$,   we set 
%
%
\begin{eqnarray}
\ell[f;g;h]({\bf a}; {\bf b}; {\bf c}) &=& \ell(f_1,f_2,f_3,\ldots; g_1,g_2,g_3,\ldots;h_1,h_2,h_3,\ldots)\label{eq:crochetabc}
\end{eqnarray}
and for series $\alpha=\alpha(x,y)$, $\beta=\beta(x,y)$ and $\gamma=\gamma(x,y)$, we also set
%
%
\begin{eqnarray}
\ell[\alpha;\beta;\gamma](x,y) &=& \ell\left(\alpha(x,y),\alpha(x^2,y^2),\alpha(x^3,y^3),\ldots; \beta(x,y),\beta(x^2,y^2),\beta(x^3,y^3),\ldots;\right. \phantom{M}\nonumber\\
&&\left.\phantom{M}\gamma(x,y),\gamma(x^2,y^2),\gamma(x^3,y^3),\ldots\right).
\label{eq:crochetxy} \end{eqnarray} 
%

%
\begin{thm} [\cite{GLL3,Timothyunlabelled}] \label{thm:WGflN} Let $\G$ be a species of graphs and $\N$ be a symmetric species of networks.
Then the edge index series and the tilde series of the species $\G\uparrow\N$ are given by
\begin{eqnarray}
W_{\G \uparrow \N}({\bf a}; {\bf b}; {\bf c}) &=&  W_\G(a_1,a_2,\ldots; W^+_\N,W^+_{\N,2},\ldots; W^-_\N,W^-_{\N,2},\ldots)\nonumber\\
&=&  W_\G[a_1;W^+_\N;W^-_\N]
\end{eqnarray}
and
\begin{eqnarray}
(\G \uparrow \N)\tilde{\ }(x,y) &=& W_\G(x,x^2,\ldots; \widetilde{\N}(x,y),\widetilde{\N}(x^2,y^2),\ldots; \widetilde{\N}_\tau(x,y),\widetilde{\N}_\tau(x^2,y^2),\ldots)\nonumber\\
&=&W_\G[x;\widetilde{\N}(x,y);\widetilde{\N}_\tau(x,y)].
\end{eqnarray}
\end{thm}
\proof  See \cite{GLL3}. \hfill\rule{2mm}{2mm}

Similarly, for a composition of networks $\M \uparrow \N$, we have the following.
\begin{thm} [\cite{Timothyunlabelled}]\label{thm:WMflN} Let $\M$ be a species of networks and $\N$ be a symmetric species of networks. 
Then the edge index series and the tilde series of the species $\M\uparrow\N$ are given by
\begin{eqnarray}
W^+_{\M \uparrow \N}({\bf a}; {\bf b}; {\bf c}) & = & 
W^+_\M[a_1;W^+_\N;W^-_\N],     \label{eq:W+MflN}\\
W^-_{\M \uparrow \N}({\bf a}; {\bf b}; {\bf c}) 
& = & 
W^-_\M[a_1;W^+_\N;W^-_\N],\label{eq:W-MflN}
\end{eqnarray}
and
\begin{eqnarray}
(\M \uparrow \N)\tilde{\ }(x,y) &=& 
W^+_\M[x;\widetilde{\N}(x,y);\widetilde{\N}_\tau(x,y)], \label{eq:MflNtilde}\\
(\M \uparrow \N)\tilde{_{\tau}}(x,y) &=& 
W^-_\M[x;\widetilde{\N}(x,y);\widetilde{\N}_\tau(x,y)].  \label{eq:MflNtautilde}
\end{eqnarray}
\end{thm}
\proof The proof is similar to that of Theorem \ref{thm:WGflN}, given in \cite{GLL3}. For (\ref{eq:W+MflN}) and (\ref{eq:W-MflN}), one uses the fact that the poles of an $\M \uparrow \N$-structure are preserved (resp. exchanged) if and only if the poles of its core are preserved (resp. exchanged). Notice that (\ref{eq:MflNtilde}) and (\ref{eq:MflNtautilde}) are consequences of (\ref{eq:W+MflN}) and (\ref{eq:W-MflN}) by virtue of Proposition \ref{propos:specDeW}.
\hfill\rule{2mm}{2mm}
\begin{propos}[\cite{GLL3,Timothyunlabelled}] Let $\B$ be a species of $2$-connected graphs, with $K_2\in\B$. Then the edge index series of the associated species of networks $\B_{0,1}$ and  $\N_\B = (1+y)\B_{0,1} - \nn$ are given by
\begin{eqnarray}
W^+_{\B_{0,1}}({\bf a}; {\bf b}; {\bf c})&=&\frac{2}{a_1^2}\frac{\partial}{\partial b_1}W_{\B}({\bf a}; {\bf b}; {\bf c}),
\label{form:18}\\
W^-_{\B_{0,1}}({\bf a}; {\bf b}; {\bf c})&=&\frac{2}{a_2}\frac{\partial}{\partial c_1}W_{\B}({\bf a}; {\bf b}; {\bf c}),
\label{form:19}\\
W^+_{\N_\B}({\bf a}; {\bf b}; {\bf c})&=&(1+b_1)W^+_{\B_{0,1}}({\bf a}; {\bf b}; {\bf c}) - 1,
\label{form:20}\\	
W^-_{\N_\B}({\bf a}; {\bf b}; {\bf c})&=&(1+c_1)W^-_{\B_{0,1}}({\bf a}; {\bf b}; {\bf c}) - 1.
\label{form:21}
\end{eqnarray}
\end{propos}
Note that for the operator $\N \mapsto y\N$, where $\N$ is a species of networks with non-adjacent poles, which consists in adding the edge $01$ to all networks in $\N$, we have
\begin{equation}
W^+_{y\N} = b_1W^+_{\N}\ \ \ \mathrm{and}\ \ \ W^-_{y\N} = c_1W^-_{\N}.
\end{equation}

For the series 
composition of networks, we have the following edge index series identities.
\begin{thm} \label{theo:WSeriesComp}
Let $\M$ and $\N$ be species of non-trivial networks. 
Then we have
\begin{eqnarray}
W^+_{\M\cdot_s\N}({\bf a}; {\bf b}; {\bf c}) & = & a_1W^+_\M({\bf a}; {\bf b}; {\bf c})W^+_\N({\bf a}; {\bf b}; {\bf c}), \label{W+MsN}
\\
W^-_{\M\cdot_s\M}({\bf a}; {\bf b}; {\bf c}) & = & a_1W^+_{\M,2}({\bf a}; {\bf b}; {\bf c}), \label{W-MsM}
\\
W^-_{\M\cdot_s\N\cdot_s\M}({\bf a}; {\bf b}; {\bf c}) & = & a_2W^+_{\M,2}({\bf a}; {\bf b}; {\bf c})W^-_\N({\bf a}; {\bf b}; {\bf c}). \label{W-MsNsM}
\end{eqnarray}
\end{thm}
\proof We can use the representations (\ref{form:35}) for the edge index series $W^+$ and $W^-$ which interpret these series as exponential generating functions of labelled structures. Thus in the first case we are interested in the exponential generating function $(\M\cdot_s\N)^+_w(x)$. Now a $\M\cdot_s\N$-structure is a pair $(M\cdot_sN, \sigma)$, where $M\cdot_sN$ is a series composition with $M\in\M$ and $N\in\N$ and $\sigma$ is a pole-preserving automorphism of $M\cdot_sN$. It is clear that the connecting vertex of the series composition is left fixed by $\sigma$ and that $\sigma$ can be restricted to pole-preserving automorphisms $\sigma_M$ and $\sigma_N$ of $M$ and $N$, respectively.  Moreover, for the weight $w$ defined by (\ref{form:weightN+}), we have
$$
w(M\cdot_sN, \sigma) = a_1w(M,\sigma_M)w(N,\sigma_N).
$$
Hence the generating functions satisfy
$$
(\M\cdot_s\N)^+_w(x) = xa_1M^+_w(x)N^+_w(x)
$$
and (\ref{W+MsN}) follows.

\begin{figure}[h] 
\begin{center} \includegraphics[height=1.30in]{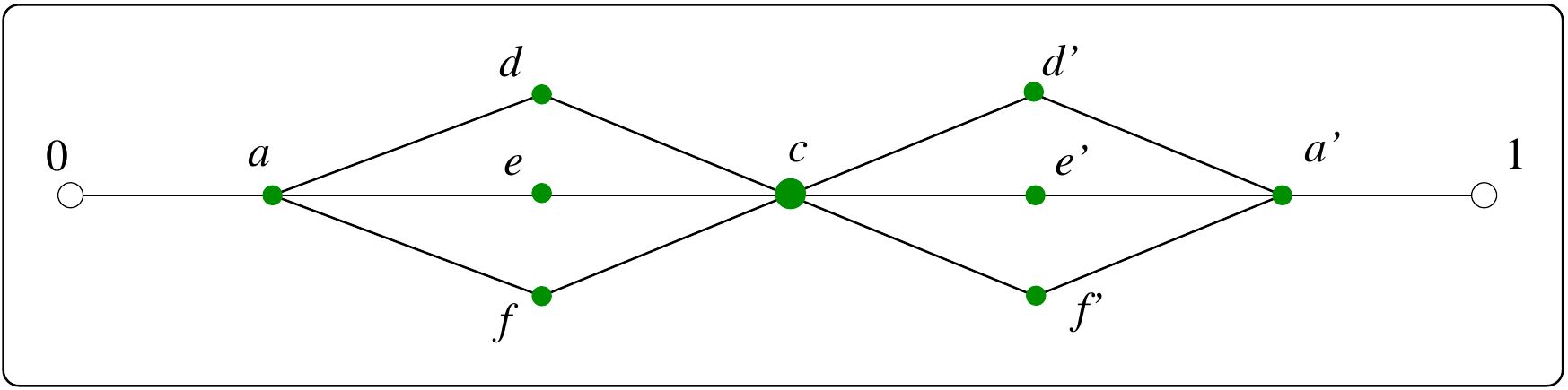}
\end{center}
\vspace{-3mm}
\caption{$\tau$-symmetric series composition of networks \label{fig:seriescomp-}}
\end{figure}

In order to prove (\ref{W-MsM}), one should enumerate structures of the form 
$(M\cdot_s M', \sigma)$, where $M\cdot_s M'$ is a $\tau$-symmetric series-composition network, with $M$ and $M'$ in $\M$, and $\sigma$ is a pole-reversing graph automorphism. In this case $\sigma$ will leave the connecting vertex $c$ fixed and will induce two network isomorphisms 
\begin{equation}
\varphi=\sigma|_M:M\tilde{\rightarrow}\tau M'\ \  \mathrm{and}\ \  \rho=\sigma^2|_M:M\tilde{\rightarrow}M.
\end{equation}
See Figure \ref{fig:seriescomp-} where the isomorphism $\varphi$ is represented as $x\mapsto x'$, for $x= a,d,e,f$. Conversely, the data of $\rho$ and $\varphi$ determines $\sigma$ since $\sigma=\varphi\cup(c)\cup\rho\circ\varphi^{-1}$. Moreover, all the vertex- and edge-cycles of $\rho$ have their lengths doubled in $\sigma$. For example, taking $\rho=(a)(d,e,f)$ in Figure~\ref{fig:seriescomp-}, we find that $\sigma=(a,a')(d,d',e,e',f,f')$.  It follows that 
$$
(\M\cdot_s\M)^-_w(x) = xa_1M^+_{(w)_2}(x^2),
$$
where we set $(w)_2(G,\sigma)=w(G,\sigma)_2$, corresponding to the plethystic notation (\ref{eq:plethisme}).

  In the case of (\ref{W-MsNsM}), the reasoning is similar.  Here a pole-reversing automorphism $\sigma$ of a series composition $M\cdot_s N\cdot_s M'$ will exchange $M$ and $M'$ and, furthermore, induce a pole-reversing automorphism of $N$ and interchange the two connecting vertices. Details are left to the reader.
\hfill\rule{2mm}{2mm}
\begin{propos} \label{prop:s-networksBis}
Let $\D=\D_\F$ be the class of non-trivial networks all of whose 3-connected components are in a given species $\F$ and let $\S$ denote the class of $s$-networks in $\D$. Then we have
\begin{equation} \label{eq:W+S}
W^+_{\S}({\bf a}; {\bf b}; {\bf c})  =  \frac{a_1 (W^+_\D)^2}{1+a_1W^+_\D}
\end{equation}
and
\begin{equation} \label{eq:W-S}
W^-_{\S}({\bf a}; {\bf b}; {\bf c}) =  \frac{(a_1+a_2W^-_\D)W^+_{\D,2}}{1+a_2W^+_{\D,2}}.
\end{equation}
\end{propos}
\proof 
From Proposition \ref{prop:s-networks}, we have 
\begin{equation} \label{eq:SDSD}
\S=(\D-\S)\cdot_s\D=\D\cdot_s\D - \S\cdot_s\D
\end{equation}
and, by Theorem \ref{theo:WSeriesComp}, 
\begin{equation}
W^+_{\S} = a_1(W^+_\D)^2 - a_1 W^+_\S W^+_\D. 
\end{equation}
Solving for $W^+_{\S}$, we obtain (\ref{eq:W+S}).  
However, formula (\ref{eq:SDSD}) can not be used for computing the edge index series $W^-_\S$ since the decomposition is not preserved by a pole-reversing automorphism.  One should rather use the more symmetric canonical decomposition (\ref{eq:SDmoinsS}) and then apply (\ref{W-MsM}) and (\ref{W-MsNsM}).  Regrouping the even and the odd $\cdot_s$-powers yields
\begin{equation}
W^-_\S  =  \frac{a_1(W^+_{\D,2}-W^+_{\S,2})}{1-a_2(W^+_{\D,2}-W^+_{\S,2})} + \frac{a_2(W^+_{\D,2}-W^+_{\S,2})(W^-_{\D}-W^-_{\S})}{1-a_2 (W^+_{\D,2}-W^+_{\S,2})},
\end{equation}
and, after simplification,
\begin{equation}
W^-_\S  = (W^+_{\D,2}-W^+_{\S,2})(a_1 + a_2W^-_\R).
\end{equation}
Formula (\ref{eq:W+S}) can then be used and the result follows.
\hfill\rule{2mm}{2mm}

\begin{thm}
Let $\N$ be a species of non-trivial networks having non-adjacent poles. 
Then the edge index series of the species of parallel compositions $E(\N)$ are given by
\begin{equation} \label{eq:W+EN}
W^+_{E(\N)}({\bf a}; {\bf b}; {\bf c})  = 
\exp\left(\sum_{m=1}^{\infty}\frac{W^+_{\N,m}}{m}\right)
\end{equation}
and 
\begin{equation} \label{eq:W-EN}
W^-_{E(\N)}({\bf a}; {\bf b}; {\bf c}) = \exp\left(\sum_{m\ \mathrm{even}}\frac{W^+_{\N,m}}{m}+\sum_{m\ \mathrm{odd}}\frac{W^-_{\N,m}}{m}\right).
\end{equation}
We also have
\begin{equation} \label{eq:W+-E2N}
W^+_{E_2(\N)}
= \frac{1}{2}((W^+_{\N})^2 + W^+_{\N,2})\ \ \mathrm{and}\ \ 
W^-_{E_2(\N)}
= \frac{1}{2}((W^-_{\N})^2 + W^+_{\N,2}).
\end{equation}
\end{thm}
\proof We  use again the representations (\ref{form:35}) for the edge index series $W^+$ and $W^-$. In the first case we are interested in the exponential generating function $E(\N)^+_w(x)$. An $E(\N)^+$-structure consists of a parallel composition of networks in $\N$ together with a network automorphism $\sigma$. This $\sigma$ induces a permutation $\sigma_0$ on the set of individual networks involved in the parallel composition. Decomposing $\sigma_0$ into (oriented) cycles yields a natural notion of \emph{connected} $E(\N)^+$-structure, namely when $\sigma_0$ is a circular permutation leading to an oriented cycle of network isomorphisms. These are known as \emph{cylindrical $m$-wreaths of networks} (see \cite{Joyal,GLL3}), 
\begin{equation}\label{eq:cm}
c_m : 
N_1 \stackrel{\varphi_1}{\longrightarrow} N_2 \stackrel{\varphi_2}{\longrightarrow}\ldots\stackrel{\varphi_{m-2}}{\longrightarrow} N_{m-1} \stackrel{\varphi_{m-1}}{\longrightarrow} N_m \stackrel{\varphi_m}{\longrightarrow} N_1,
\end{equation} 
where $m\geq1$.
Let $\K_m(\N)$ denote the species of cylindrical $m$-wreaths of networks and $\K_m^{\bullet}(\N)$ species of \emph{rooted} cylindrical $m$-wreaths of networks, where one network is distinguished from the others.
In fact the description (\ref{eq:cm}) includes a rooting at $N_1$. In the unrooted case, all the possible rootings are considered equivalent. It follows that any $E(\N)^+$-structure can be seen as an assembly of (unrooted) cylindrical wreaths of networks and we have a weighted species isomorphism
\begin{equation} \label{eq:ENKmN}
E(\N)^+_w = E\left(\sum_{m\geq1} \K_m(\N)_w\right).
\end{equation} 
Given a rooted cylindrical $m$-wreath of networks $c_m$ in $\K_m^{\bullet}(\N)$, of the form (\ref{eq:cm}), the composite 
$\varphi_0 = \varphi_m \circ \varphi_{m-1} \circ \ldots \varphi_2 \circ \varphi_1$ is an automorphism of $N_1$, and we obtain a $\N^+$-structure $(N_1,\varphi_0)$. Moreover the sequence of network isomorphisms $(\varphi_1,\ldots,\varphi_{m-1})$ can be encoded in a set of lists of length $m$, $(u_1, u_2, \dots, u_m)$, where $u_1$ runs over the underlying set of $N_1$ and $u_{i+1}= \varphi_i(u_i)$, $i=1\ldots m-1$, and we can consider the $\N^+$-structure $(N_1,\varphi_0)$ to ``live" on this set of lists. In other words, what we have obtained is an $\N^+(X^m)$-structure. Since the isomorphism $\varphi_m$ can be recovered from $\varphi_0$ and the other isomorphisms $\varphi_i$, this correspondence is bijective. Moreover the weight of the connected $E(\N)^+$-structure $c_m$ is given by $(w)_m(N_1,\varphi_0)$ since all the cycle lengths of $\varphi_0$ are multiplied by $m$ in $\varphi_0$. Hence we have an isomorphism of weighted species (see also Proposition 14 of \cite{GLL3})
\begin{equation} \label{eq:KmNbullet}
\K_m^{\bullet}(\N)_{w} = \N_{(w)_m}^+(X^m)
\end{equation}
and the exponential generating function equality
\begin{equation} \label{eq:KmNbulletx}
\K_m(\N)_{w}(x) = \frac{1}{m}\K_m^{\bullet}(\N)_{w}(x) = \frac{1}{m}\N_{(w)_m}^+(x^m).
\end{equation}
Using (\ref{eq:ENKmN}) and the classical exponential formula, we find that
\begin{equation} \label{eq:ENKmNx}
E(\N)^+_w(x) = \exp\left(\sum_{m\geq1} \frac{1}{m}\N_{(w)_m}^+(x^m)\right) 
\end{equation} 
and (\ref{eq:W+EN}) follows.

For (\ref{eq:W-EN}), one should compute  the exponential generating function $E(\N)^-_w(x)$. An $E(\N)^-$-structure consists of a parallel composition of networks in $\N$ together with a pole-reversing automorphism $\sigma$. Here two kinds of connected components can occur.  The first kind arises from a cylindrical $m$-wreath of networks such as (\ref{eq:cm}), with $m$ even, which is reinterpreted as a sequence of pole-reversing network isomorphisms
\begin{equation}\label{eq:cm-pair}
N_1 \stackrel{\varphi_1}{\longrightarrow} \tau N_2 \stackrel{\varphi_2}{\longrightarrow} N_3\stackrel{\varphi_3}{\longrightarrow}\ldots\stackrel{\varphi_{m-2}}{\longrightarrow} N_{m-1} \stackrel{\varphi_{m-1}}{\longrightarrow} \tau N_m \stackrel{\varphi_m}{\longrightarrow} N_1.
\end{equation}
Here also the composite $\varphi_0 = \varphi_m \circ \varphi_{m-1} \circ \ldots \varphi_2 \circ \varphi_1$ is an automorphism of $N_1$ and this accounts for the fisrt term on the right-hand side of (\ref{eq:W-EN}). 

  The second kind of connected component corresponds to a \emph{M\"obius $m$-wreath of networks}, with $m$ odd, which is defined as a sequence of network isomorphisms 
$N_1 \stackrel{\varphi_1}{\longrightarrow} N_2 \stackrel{\varphi_2}{\longrightarrow}\ldots\stackrel{\varphi_{m-2}}{\longrightarrow} N_{m-1} \stackrel{\varphi_{m-1}}{\longrightarrow} N_m$ 
followed by a pole-reversing isomorphism $\varphi_m : N_m \longrightarrow \tau N_1$, which can be reinterpreted as a sequence of pole-reversing  isomorphisms
\begin{equation}\label{eq:cm-pair}
N_1 \stackrel{\varphi_1}{\longrightarrow} \tau N_2 \stackrel{\varphi_2}{\longrightarrow} N_3\stackrel{\varphi_3}{\longrightarrow}\ldots\stackrel{\varphi_{m-2}}{\longrightarrow}\tau N_{m-1} \stackrel{\varphi_{m-1}}{\longrightarrow} N_m \stackrel{\varphi_m}{\longrightarrow} \tau N_1.
\end{equation}
Notice that the composite
$\varphi_0 = \varphi_m \circ \varphi_{m-1} \circ \ldots \varphi_2 \circ \varphi_1$ is a pole-reversing automorphism of $N_1$ and this accounts for the second term on the right-hand side of (\ref{eq:W-EN}).  See also \cite{GLL3}.

  The proof of (\ref{eq:W+-E2N}) relies on the fact that 
\begin{equation}
E_2(\N)^+_w = E_2(\N^+_w) + \K_2(\N)_w \ \  \mathrm{and} \ \  E_2(\N)^-_w = E_2(\N^-_w) + \K_2(\N)_w\,.
\end{equation}
Details are left to the reader. 
\hfill\rule{2mm}{2mm}

\smallskip\noindent
\textbf{Remark.} A proof of (\ref{eq:W+EN}) and (\ref{eq:W-EN}) obtained by expressing a parallel composition as an $\uparrow$-composition whose core is a ``network'' with parallel edges and no internal vertices appears in \cite{Timothyunlabelled}.
%
%
\section{Enumerative applications} 
Again let $\F$ be a given class of 3-connected graphs and let $\B = \B_{\F}$ denote the class of 2-connected graphs all of whose 3-connected components are in $\F$. 
Also let $\D = \D_{\F}$ denote the class of networks all of whose 3-connected components are in $\F$. Given $\F$, the species $\D$ can be determined recursively, as well as its associated series, from the fundamental relations of Theorem \ref{theo:reseaux}. Using the dissymmetry theorem, the species $\B$ and its series can also be determined. 

The formulas of the previous sections can be applied to obtain both the labelled and unlabelled enumeration of species of 2- or 3-connected graphs.  
The labelled enumeration is usually simpler since it is not necessary to use the Dissymmetry Theorem in order to unroot the structures and since also the composition formulas are simpler for the exponential generating functions.
For unlabelled enumeration, the formulas are more delicate and are reviewed below. Some standard applications and some new ones are also presented. 
\subsubsection{Labelled enumeration}  For the network species $\D=\D_\F$, 
we deduce the following functional equation for the exponential generating function: 
\begin{equation} \label{eq:RxdeRx}
\D(x,y) = (1+y)\exp\left(\F_{0,1}(x,\D(x,y)) + \frac{x\D^2(x,y)}{1+x\D(x,y)}\right)-1.
\end{equation}
Setting $\g(x,y) = \F_{0,1}(x,y) + \frac{xy^2}{1+xy}$, and $\z(x,y)= \D(x,y)^{<-1>_y}$, the compositional inverse of $\D(x,y)$ with respect to $y$, 
we have, from (\ref{eq:RxdeRx}),
\begin{equation}
1+\D(x,y) = (1+y)\exp(\g(x,\D(x,y)))
\end{equation}
and
\begin{equation}\label{eq:zetaxy}
1+\z(x,y) = (1+y)\exp(-\g(x,y)).
\end{equation}
Notice that $\z(x,y)$ is of the form
\begin{equation}
\z(x,y)=y\left(1 + (1+y) \frac{\exp(-g(x,y))-1}{y}\right)
\end{equation}
so that Lagrange inversion can be used to find $\D(x,y)$, knowing $\F_{0,1}(x,y)$.

Conversely, taking logarithms in (\ref{eq:zetaxy}) yields $\F_{0,1}(x,y)$ in terms of $\D(x,y)$:
\begin{equation}
\F_{0,1}(x,y)=\log\frac{1+y}{1+\z(x,y)}-\frac{xy^2}{1+xy}.
\end{equation}
Finally, note that 
\begin{equation}
\F(x,y)= \frac{x^2}{2}\int \F_{0,1}(x,y)dy \ \ \mathrm{and}\ \  
\B(x,y)=\frac{x^2}{2}\int\frac{1+\D(x,y)}{1+y}dy. \label{eq:BDxy}
\end{equation}

This is essentially the approach used in \cite{Timothy} for the enumeration of labelled 3-connected graphs, starting with 1- and 2-connected graphs, and in \cite{Bender}, where labelled 2-connected planar graphs are enumerated, starting from 3-connected planar graphs. 
%
%
\subsubsection{Unlabelled enumeration} \label{sec:unlabelled}
We introduce the following abbreviations for the
edge index series of the species $\R$, $\S$, of $s$-networks,  and $\H=\F_{0,1}\uparrow\R$, of $h$-networks:
\begin{equation}
\rho^+({\bf a}; {\bf b}; {\bf c})= W^+_\R({\bf a}; {\bf b}; {\bf c}), \ \ \ 
\rho^-({\bf a}; {\bf b}; {\bf c})= W^-_\R({\bf a}; {\bf b}; {\bf c}),
\end{equation}
\begin{equation}
\sigma^+({\bf a}; {\bf b}; {\bf c})= W^+_S({\bf a}; {\bf b}; {\bf c}), \ \ \ 
\sigma^-({\bf a}; {\bf b}; {\bf c})= W^-_S({\bf a}; {\bf b}; {\bf c}),
\end{equation}
and
\begin{equation}
\eta^+({\bf a}; {\bf b}; {\bf c})= W^+_H({\bf a}; {\bf b}; {\bf c}), \ \ \ 
\eta^-({\bf a}; {\bf b}; {\bf c})= W^-_H({\bf a}; {\bf b}; {\bf c}).
\end{equation}
It follows from equations (\ref{eq:SXD2D}),  (\ref{eq:HF01D}) and (\ref{eq:RdeR}) and from the properties of the edge index series, that  
\begin{eqnarray}
\rho^+ & = &(1+b_1)\exp\left(\sum_{i=1}^{\infty}\frac{\eta^+_i+\sigma^+_i}{i}\right)-1,\label{eq:rhoderho+}\\
\rho^- & = & (1+c_1)\exp\left(\sum_{i\ \mathrm{even}}\frac{\eta^+_i+\sigma^+_i}{i}+\sum_{i\ \mathrm{odd}}\frac{\eta^-_i+\sigma^-_i}{i}
\right)-1, \label{eq:rhoderho-}\\
\eta^+ &=& W^+_{\F_{0,1}}[a_1;\rho^+;\rho^-],\ \ \ 
\eta^- =\ W^-_{\F_{0,1}}[a_1;\rho^+;\rho^-], \label{eq:eta}\\
\sigma^+ &=& \frac{a_1(\rho^+)^2}{1+a_1\rho^+}, \ \ \ \ \ \ \ 
\sigma^- =\ \frac{(a_1+a_2\rho^-)\rho^+_2}{1+a_2\rho^+_2}. \label{eq:sigma}
\end{eqnarray}
These equations make it possible to compute recursively the series $\rho^+$, $\rho^-$, $\eta^+$, $\eta^-$, $\sigma^+$ and $\sigma^-$, knowing $W^+_{\F_{0,1}}$ and $W^-_{\F_{0,1}}$. The dissymmetry formula (\ref{eq:dissymBis}) will then yield $W_\B$, knowing $W_\F\,$:
\begin{eqnarray}
W_\B({\bf a}; {\bf b}; {\bf c}) & = & W_\F[a_1;\rho^+;\rho^-]\ +\  W_\C[a_1;\rho^+-\sigma^+;\rho^--\sigma^-] \nonumber  \\ 
& &  +\ \frac{a_1^2}{2} \left(\rho^+ -(\eta^++\sigma^+)(\rho^++1) - \frac{1}{2}(\sigma^+_2 - (\sigma^+)^2) \right) \nonumber  \\
& & +\ \frac{a_2}{2} \left(\rho^- -(\eta^-+\sigma^-)(\rho^-+1) - \frac{1}{2}(\sigma^+_2 - (\sigma^-)^2) \right), \label{eq:WdissymBis}
\end{eqnarray}
where $W_\C$ is given by (\ref{eq:WC}).
If only the tilde generating functions are desired, for the unlabelled enumeration, equations (\ref{eq:rhoderho+} -- \ref{eq:sigma}) yield the following:
\begin{eqnarray}
\widetilde{\R}(x,y) & = &(1+y)\exp\left(\sum_{i=1}^{\infty}\frac{(\H+\S)\tilde{\ }(x^i,y^i)}{i}\right)-1,\label{eq:Rtildexy}\\
\widetilde{\R}_\tau(x,y) & = & (1+y)\exp\left(\sum_{i\ \mathrm{even}}\frac{(\H+\S)\tilde{\ }(x^i,y^i)}{i}+\sum_{i\ \mathrm{odd}}\frac{(\H+\S)\tilde{_\tau }(x^i,y^i)}{i}
\right)-1, \label{eq:Rtautildexy} 
\end{eqnarray}
\begin{eqnarray}
\widetilde{\H}(x,y) &=& W^+_{\F_{0,1}}[x;\widetilde{\R}(x,y);\widetilde{\R}_\tau(x,y)],\label{eq:Htildexy}\\
\widetilde{\H}_\tau(x,y) &=&\ W^-_{\F_{0,1}}[x;\widetilde{\R}(x,y);\widetilde{\R}_\tau(x,y)], \label{eq:Htautildexy}
\end{eqnarray}
where the notation of (\ref{eq:crochetxy}) is used, and 
\begin{eqnarray}
\widetilde{\S}(x,y) &=& \frac{x\widetilde{\R}^2(x,y)}{1+x\widetilde{\R}(x,y)}, \label{eq:Stildexy}\\
\widetilde{\S}_\tau(x,y)&=& \frac{(x+x^2\widetilde{\R}_\tau(x,y))\widetilde{\R}(x^2,y^2)}{1+x^2\widetilde{\R}(x^2,y^2)}. \label{eq:Stautildexy}
\end{eqnarray}
Finally, equation (\ref{eq:WdissymBis}) 
gives the following dissymmetry formula:
\begin{eqnarray}
\widetilde{\B}(x,y) & = & W_\F[x;\;\widetilde{\R};\;\widetilde{\R}_\tau]\ +\  W_\C[x;\;\widetilde{\R}-\widetilde{\S};\;\widetilde{\R}_\tau-\widetilde{\S}_\tau]\nonumber  \\ 
& &  +\ \frac{x^2}{2} \left(\widetilde{\R} -(\widetilde{\H}+\widetilde{\S})(\widetilde{\R}+1) - \widetilde{\S}_2 + \frac{1}{2}\widetilde{\S}^{\,2} \right) \nonumber  \\
& & +\ \frac{x^2}{2} \left(\widetilde{\R}_\tau -(\widetilde{\H}_\tau+\widetilde{\S}_\tau)(\widetilde{\R}_\tau+1) + \frac{1}{2} \widetilde{\S}_\tau^{\,2} \right), \label{eq:dissymBisxy}
\end{eqnarray}
where $\widetilde{\S}_2(x,y) = \widetilde{\S}(x^2,y^2)$.

%
%
%
\subsection{3-connected graphs} 
The first application of the above formulas is the enumeration of unlabelled 3-connected graphs in 1982 (see \cite{Timothyunlabelled}, \cite{RobinsonWalsh}).
In this case, where $\F=\F_a$ and $\B=\B_a$, the species of all 3-connected and all 2-connected (simple) graphs, respectively, it is possible to compute the edge index series $W_\B$ directly, going from all graphs, to connected graphs and then to 2-connected graphs. Since $\D=(1+y)\B_{0,1}-1$, the edge index series $\rho^+=W^+_\R$ and $\rho^-=W^-_\R$ can be readily computed, as well as $\sigma^+=W^+_\S$ and $\sigma^-=W^-_\S$, using (\ref{eq:sigma}), and then $\eta^+=W^+_\H$ and $\eta^-=W^-_\H$, recursively, using equations (\ref{eq:rhoderho+}) and (\ref{eq:rhoderho-}). The edge index series 
$W_\F$ is then extracted recursively, using 
the dissymmetry formula (\ref{eq:WdissymBis}) and the generating function $\widetilde{\F}(x,y)$ is then immediately deduced.
In \cite{RobinsonWalsh}, the computations are greatly simplified by introducing
two auxiliary series $\beta(x,y)$ and $\gamma(x,y)$ satisfying
\begin{equation}
\rho^+[x,\beta(x,y),\gamma(x,y)] = y\, \ \ \ \ \rho^-[x,\beta(x,y),\gamma(x,y)] = y. 
\end{equation}
The readers are referred to \cite{RobinsonWalsh} for more details.  See \cite{Timothyunlabelled} and \cite{RobinsonWalsh} for tables.
%
%
\subsection{Series-parallel graphs and networks} 
At the other extreme lies the case where $\F=0$, where the corresponding species  of 2-connected graphs is the class $\B=\SPG$ of \emph{series-parallel graphs}. Thus a series-parallel graph $G$ is a 2-connected graph all of whose 3-components are polygons.  $G$ can also be characterized by the fact it contains no subdivision of $K_4$.   

The corresponding species of networks is the class $\D=\DSP$ of \emph{series-parallel networks}. An example of a series-parallel network is given in Example \ref{fig:reseauxexemple}. This class $\D$ is defined recursively by the functional equation 
\begin{equation} \label{eq:RdeRSP}
\D = (1+y)E(\frac{X\D^2}{1+X\D})- \mathrm{\nn},
\end{equation}
which is a specialization to $\F=0$, of equation (\ref{eq:RdeR}).
We set $\rho^+ = W^+_{\R}({\bf a}; {\bf b}; {\bf c})$ and $\rho^- = W^-_{\R}({\bf a}; {\bf b}; {\bf c})$. Equations (\ref{eq:rhoderho+}) and (\ref{eq:rhoderho-}) then imply the following:

\begin{corol} For the edge index series $\rho^+$ and $\rho^-$ of the species $\D=\DSP$ of series-parallel networks, we have the system of equations 
\begin{eqnarray}
\rho^+ & = &(1+b_1) \exp\left(\sum_{i=1}^{\infty}\frac{1}{i}\frac{a_i(\rho^+_i)^2}{1+a_i\rho^+_i}\right) - 1, \label{eq:rhoderho+SP}\\
\rho^- & = &(1+c_1) \exp\left(\sum_{i\ \mathrm{even}} \frac{1}{i} \frac{a_i(\rho^+_i)^2}{1+a_i\rho^+_i} +\sum_{i\ \mathrm{odd}}\frac{1}{i}\frac{(a_i+a_{2i}\rho^-_i)\rho^+_{2i}}{1+a_{2i}\rho^+_{2i}}\right) - 1.  \label{eq:rhoderho-SP}
\end{eqnarray}
\end{corol}
These edge index series can be computed recursively and the edge index series $W_{\SPG}$ of series parallel graphs can be deduced by specializing equation (\ref{eq:WdissymBis}):
\begin{eqnarray}
W_{\SPG}({\bf a}; {\bf b}; {\bf c}) & = &  W_\C[a_1;\rho^+-\sigma^+;\rho^--\sigma^-] 
+\ \frac{a_1^2}{2} (\rho^+ -\sigma^+ -\sigma^+\rho^+ - \frac{1}{2}(\sigma^+_2 - (\sigma^+)^2) ) \nonumber  \\
& & +\ \frac{a_2}{2} (\rho^- -\sigma^- -\sigma^-\rho^- - \frac{1}{2}(\sigma^+_2 - (\sigma^-)^2)), \label{eq:WdissymBisSP}
\end{eqnarray}
where $\sigma^+$ and $\sigma^-$  are defined by (\ref{eq:sigma}). The exponential and tilde generating functions are then immediately obtained.  For example, for unlabelled series-parallel graphs counted according to the number of vertices (where we set $y=1$), we find 
\begin{eqnarray} 
\widetilde{\SPG}(x) &=& x^2 + x^3 + 2x^4 + 5x^5 + 15x^6 + 51x^7 + 230x^8 + 1142x^9 \nonumber\\
&& +\: 6369x^{10} + 37601x^{11} + 232259x^{12} + 1476120x^{13} + 9599522x^{14} + \cdots \phantom{MM}
\end{eqnarray} 
Similarly, for unlabelled series-parallel networks counted according to the number of internal vertices, we find that
\begin{eqnarray} 
\widetilde{\DSP}(x) &=&
1+2\,x+8\,{x}^{2}+38\,{x}^{3}+208\,{x}^{4}+1220\,{x}^{5}+7592\,{x}^{6
}+49006\,{x}^{7}\nonumber \\ && +325686\,{x}^{8} 
 + 2212112\,{x}^{9}+15290182\,{x}^{10}+
107191458\,{x}^{11}+760349722\,{x}^{12}\nonumber \\ && +5447100396\,{x}^{13}+
39354320204\,{x}^{14}+
\cdots
\end{eqnarray} 
and for those that are $\tau$-symmetric, we have
\begin{eqnarray} 
\widetilde{\DSP}_\tau(x) &=&
1+2\,x+4\,{x}^{2}+10\,{x}^{3}+24\,{x}^{4}+64\,{x}^{5}+168\,{x}^{6}+
458\,{x}^{7}+1250\,{x}^{8}+3492\,{x}^{9}\nonumber \\ && +9734\,{x}^{10}+27582\,{x}^{11}+78078\,{x}^{12}+223644\,{x}^{13}+639948\,{x}^{14}+ \cdots
\end{eqnarray} 

In comparing with the existing literature on series-parallel networks, recall that here, we are considering 01-networks without parallel edges.
%
%
\subsection{2-connected planar graphs and strongly planar networks.} \label{sec:planar}
Let $\F_P$ denote the species of 3-connected planar graphs. Then the corresponding species $\B(\F_P)=\B_P$ is the class of 2-connected planar graphs and $\R(\F_P)=\N_P$ is the class of \emph{strongly planar} networks, that is of non-trivial networks $N$ such that $N\cup01$ is planar.  As before, we have
\begin{equation*}
\N_P = (1+y)(\B_{P})_{01} - \nn.
\end{equation*}
Here we compute all the edge index series and generating functions of these species up to 14 vertices.  

The enumeration of unlabelled planar graphs is one of the classical fundamental open problems in graph theory and combinatorics. 
A bridge between planar graph enumeration and planar map enumeration is provided by the fact that a 3-connected planar graph admits a unique embedding on the sphere up to homeomorphisms that either preserve or reverse the orientation of the sphere and that any graph automorphism of a 3-connected planar graph is a map automorphism of the corresponding unsensed map, and conversely (see \cite{Mani}).  
%
%
%
%
\begin{table}[!b] 
\centerline {\scriptsize
\begin{tabular}[t]{|| r | r | r || r | r | r || r | r | r ||}
\hline
$n$ & $m$ & $g_{n,m}$ & $n$ & $m$ & $g_{n,m}$ & $n$ & $m$ & $g_{n,m}$\\
\hline \hline
\hline
2 & 1 & 1 & 10 & 10 & 1 & 13 & 13 & 1\\
\hline
3 & 3 & 1 & 10 & 11 & 9 & 13 & 14 & 15\\
\hline
4 & 4 & 1 & 10 & 12 & 121 & 13 & 15 & 428\\
\hline
4 & 5 & 1 & 10 & 13 & 1018 & 13 & 16 & 8492\\
\hline
4 & 6 & 1 & 10 & 14 & 5617 & 13 & 17 & 107771\\
\hline
5 & 5 & 1 & 10 & 15 & 20515 & 13 & 18 & 903443\\
\hline
5 & 6 & 2 & 10 & 16 & 52068 & 13 & 19 & 5287675\\
\hline
5 & 7 & 3 & 10 & 17 & 94166 & 13 & 20 & 22514501\\
\hline
5 & 8 & 2 & 10 & 18 & 123357 & 13 & 21 & 71869047\\
\hline
5 & 9 & 1 & 10 & 19 & 116879 & 13 & 22 & 175632924\\
\hline
6 & 6 & 1 & 10 & 20 & 79593 & 13 & 23 & 333410770\\
\hline
6 & 7 & 3 & 10 & 21 & 37859 & 13 & 24 & 496146048\\
\hline
6 & 8 & 9 & 10 & 22 & 12066 & 13 & 25 & 581318637\\
\hline
6 & 9 & 13 & 10 & 23 & 2306 & 13 & 26 & 536073583\\
\hline
6 & 10 & 11 & 10 & 24 & 233 & 13 & 27 & 386948719\\
\hline
6 & 11 & 5 & 11 & 11 & 1 & 13 & 28 & 216020293\\
\hline
6 & 12 & 2 & 11 & 12 & 11 & 13 & 29 & 91369743\\
\hline
7 & 7 & 1 & 11 & 13 & 189 & 13 & 30 & 28288016\\
\hline
7 & 8 & 4 & 11 & 14 & 2210 & 13 & 31 & 6047730\\
\hline
7 & 9 & 20 & 11 & 15 & 16650 & 13 & 32 & 797583\\
\hline
7 & 10 & 49 & 11 & 16 & 83105 & 13 & 33 & 49566\\
\hline
7 & 11 & 77 & 11 & 17 & 289532 & 14 & 14 & 1\\
\hline
7 & 12 & 75 & 11 & 18 & 727243 & 14 & 15 & 18\\
\hline
7 & 13 & 47 & 11 & 19 & 1347335 & 14 & 16 & 616\\
\hline
7 & 14 & 16 & 11 & 20 & 1861658 & 14 & 17 & 15350\\
\hline
7 & 15 & 5 & 11 & 21 & 1926664 & 14 & 18 & 243897\\
\hline
8 & 8 & 1 & 11 & 22 & 1485235 & 14 & 19 & 2550530\\
\hline
8 & 9 & 6 & 11 & 23 & 841152 & 14 & 20 & 18598574\\
\hline
8 & 10 & 40 & 11 & 24 & 339390 & 14 & 21 & 98777626\\
\hline
8 & 11 & 158 & 11 & 25 & 92751 & 14 & 22 & 394640925\\
\hline
8 & 12 & 406 & 11 & 26 & 15362 & 14 & 23 & 1214212848\\
\hline
8 & 13 & 662 & 11 & 27 & 1249 & 14 & 24 & 2926745166\\
\hline
8 & 14 & 737 & 12 & 12 & 1 & 14 & 25 & 5594151239\\
\hline
8 & 15 & 538 & 12 & 13 & 13 & 14 & 26 & 8546948259\\
\hline
8 & 16 & 259 & 12 & 14 & 292 & 14 & 27 & 10481908901\\
\hline
8 & 17 & 72 & 12 & 15 & 4476 & 14 & 28 & 10324262525\\
\hline
8 & 18 & 14 & 12 & 16 & 44297 & 14 & 29 & 8139338353\\
\hline
9 & 9 & 1 & 12 & 17 & 290680 & 14 & 30 & 5095275794\\
\hline
9 & 10 & 7 & 12 & 18 & 1333029 & 14 & 31 & 2497740781\\
\hline
9 & 11 & 70 & 12 & 19 & 4434175 & 14 & 32 & 937658515\\
\hline
9 & 12 & 426 & 12 & 20 & 10992850 & 14 & 33 & 260094850\\
\hline
9 & 13 & 1645 & 12 & 21 & 20663187 & 14 & 34 & 50215417\\
\hline
9 & 14 & 4176 & 12 & 22 & 29764598 & 14 & 35 & 6022143\\
\hline
9 & 15 & 7307 & 12 & 23 & 32990517 & 14 & 36 & 339722\\
\hline
9 & 16 & 8871 & 12 & 24 & 28087447 & \multicolumn{3}{c}{ }\\
\cline {1-6}
9 & 17 & 7541 & 12 & 25 & 18199252 & \multicolumn{3}{c }{ }\\
\cline {1-6}
9 & 18 & 4353 & 12 & 26 & 8814281 & \multicolumn{3}{c }{ }\\
\cline {1-6}
9 & 19 & 1671 & 12 & 27 & 3088000 & \multicolumn{3}{c }{ }\\
\cline {1-6}
9 & 20 & 378 & 12 & 28 & 740272 & \multicolumn{3}{c }{ }\\
\cline {1-6}
9 & 21 & 50 & 12 & 29 & 108597 & \multicolumn{3}{c }{ }\\
\cline {1-6}
\multicolumn{3}{c ||}{ } & 12 & 30 & 7595 & \multicolumn{3}{c }{ }\\
\cline{4-6}
\end{tabular} }
	\caption{The number $g_{n,m}$ of unlabelled 2-connected planar graphs having $n$ vertices and $m$ edges.\label{table:bpnm}}
\end{table}

In order to obtain the desired series, we used the program \textit{Plantri} \cite{plantri} to generate all the 3-connected 
planar maps (alias plane graphs) up to 14 vertices.  In fact, version 4.2 of Plantri contains an option which yields the graphs together with all their automorphisms and we computed their edge cycle indices, both cylindrical and M\"obius. In this way, the edge index series $W_{\F_P}$ of 3-connected planar graphs was secured up to 14 vertices.  This yields successively  the edge index series $W^+$ and $W^-$ for the species $(\F_{P})_{01}$, $\R(\F_P)=\N_P$, $\S$ and $\H$ (recursively), and eventually the edge index series $W_{\B_P}$ of 2-connected planar graphs, using the formulas of section \ref{sec:unlabelled}. The corresponding generating functions are immediately deduced.  Thus, we extended to 14 vertices the generating function $\widetilde{\B_P}(x,y)$ of unlabelled 2-connected planar graphs and to 12 internal vertices the generating functions $\widetilde{\N_P}(x,y)$ and $\widetilde{\N_P}_\tau(x,y)$ of unlabelled strongly planar networks. The coefficients of $\widetilde{\B_P}(x,y)$ are given in Table \ref{table:bpnm}.
Setting $y=1$, we have
\begin{eqnarray} 
\widetilde{\B_P}(x) &=& x^2 + x^3 + 3x^4 + 9x^5 + 44x^6 + 294x^7 + 2893x^8 + 36496x^9 + 545808x^{10} \nonumber\\
&& +\: 9029737x^{11} + 159563559x^{12} + 2952794985x^{13} + 56589742050x^{14} + \cdots \phantom{MM}
\end{eqnarray} 
\begin{eqnarray} 
\widetilde{\N_P}(x) &=&
1+2\,x+10\,{x}^{2}+72\,{x}^{3}+696\,{x}^{4}+8530\,{x}^{5}+124926\,{x}
^{6}\nonumber\\
&& +2068888\,{x}^{7}+37204942\,{x}^{8}+708076350\,{x}^{9}+14038364914
\,{x}^{10}\phantom{MMMMM}\nonumber\\
&& +287091103062\,{x}^{11} +6016760068874\,{x}^{12} 
+ \cdots 
\end{eqnarray} 
and
\begin{eqnarray} 
\widetilde{\N_{P}}_\tau(x) \hspace{-2mm} &=& \hspace{-2mm}
1+2\,x+6\,{x}^{2}+20\,{x}^{3}+96\,{x}^{4}+470\,{x}^{5}+3074\,{x}^{6}+
23408\,{x}^{7} +243482\,{x}^{8}\nonumber\\
&&+3221018\,{x}^{9} +51729286\,{x}^{10}+ 929983374\,{x}^{11} +17911049418\,{x}^{12}
+ \cdots \phantom{MM}
\end{eqnarray} 
\textbf{Remark.}  Finding the edge index series for 3-connected planar maps without generating them is still an open problem. Wormald enumerated planar maps up to homeomorphism (see \cite{Wo1,Wo2} for 1-connected maps) without using cycle indices.
\subsection{$K_{3,3}$-free 2-connected graphs} \label{sec:K33free}
A graph is called $K_{3,3}$-\emph{free} if it contains no subdivison of $K_{3,3}$ or, equivalently, if it has no minor isomorphic to $K_{3,3}$. 
As mentioned in the introduction, a theorem of Wagner \cite{Wagner} and of Kelmans \cite{Kelmans} implies that 
if we take $\F = \F_P + K_5$, then the corresponding species $\B=\B(\F)$ of 2-connected graphs all of whose 3-connected components are planar or isomorphic to $K_5$ is the class of $K_{3,3}$-free 2-connected graphs. Computations are similar to the preceding section \ref{sec:planar}. For example we find that
\begin{eqnarray} 
\widetilde{\B}(x)\hspace{-2mm} &=& \hspace{-2mm}
{x}^{2}+{x}^{3}+3\,{x}^{4}+10\,{x}^{5}+46\,{x}^{6}+308\,{x}^{7}+2997
\,{x}^{8}+37471\,{x}^{9}+556637\,{x}^{10} \nonumber\\
&& +9171526\,{x}^{11}+161679203
\,{x}^{12}+2987857791\,{x}^{13}+57218439783\,{x}^{14}+
\cdots \phantom{MN}
\end{eqnarray} 
\begin{eqnarray} 
\widetilde{\D}(x) &=&
1+2\,x+10\,{x}^{2}+74\,{x}^{3}+718\,{x}^{4}+8786\,{x}^{5}+128006\,{x}
^{6}\nonumber\\
&& +2108610\,{x}^{7}+37767136\,{x}^{8}+716900760\,{x}^{9}+14191084858
\,{x}^{10}\phantom{MMMMM}\nonumber\\
&& +289958295858\,{x}^{11}+6074048514588\,{x}^{12} + \cdots 
\end{eqnarray} 
and
\begin{eqnarray} 
\widetilde{\D}_\tau(x) \hspace{-2mm} &=& \hspace{-2mm}
1+2\,x+6\,{x}^{2}+22\,{x}^{3}+102\,{x}^{4}+518\,{x}^{5}+3362\,{x}^{6}
+25890\,{x}^{7}+267988\,{x}^{8}\nonumber\\
&& +3524132\,{x}^{9}+56099830\,{x}^{10}+ 1001483346\,{x}^{11}+19189524860\,{x}^{12}
+ \cdots \phantom{MM}
\end{eqnarray} 
%

%
%
\subsection{$K_{3,3}$-free  projective-planar and toroidal graphs.}
Let $\PP$ denote the species of projective-planar graphs which are $K_{3,3}$-free, 2-connected and non-planar. 
In \cite{GLL} we proved the following structural characterization of $\PP$:
\begin{theorem}[\cite{GLL}]  The species $\PP$ of $K_{3,3}$-free $2$-connected non-planar projective-planar graphs can be expressed as a canonical composition
\begin{equation}
\PP = K_5\uparrow \N_P,
\end{equation} 
where $\N_P$ denotes the species of strongly planar networks. \hfill\rule{2mm}{2mm}
\end{theorem}
\begin{figure}[h]
\begin{center}
	\includegraphics[width=2.5in]{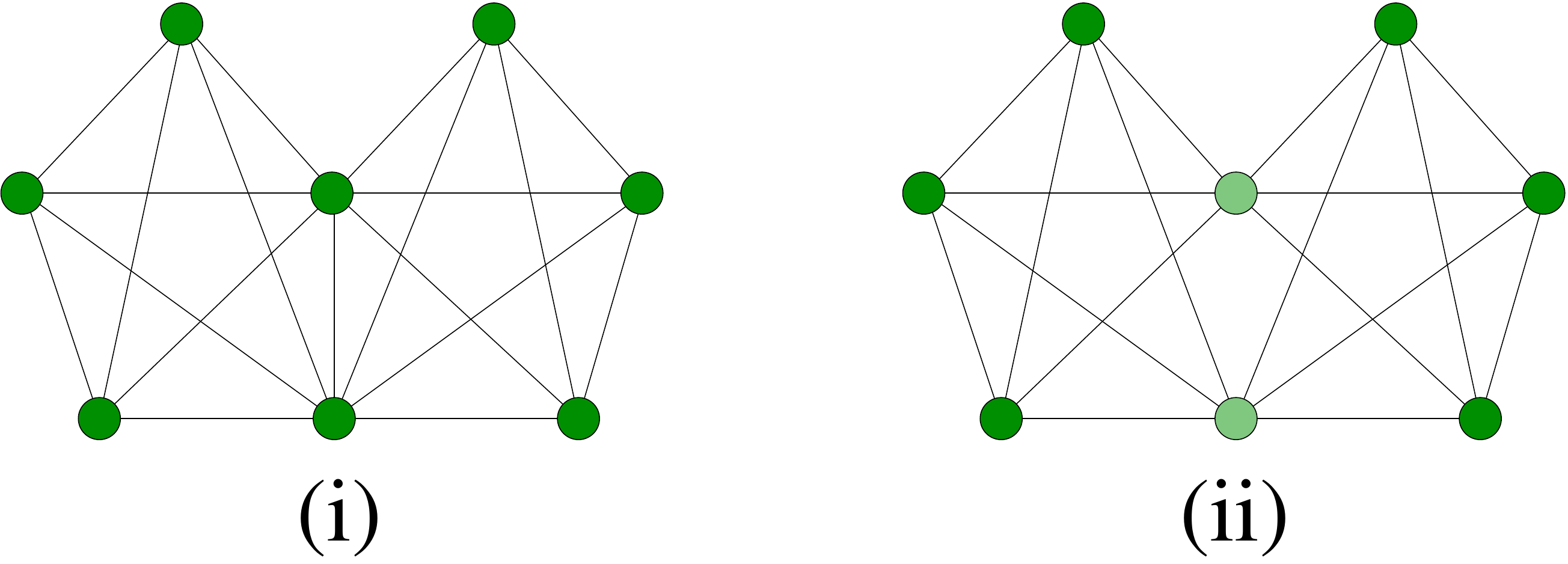}
\end{center}
\vspace{-5mm}
	\caption{(i) The graph $M$\ \ \ (ii) The graph $M^*$ \label{fig:Mgraph}}
\end{figure}
A similar characterization is provided in \cite{GLL2} for the species $\T$ of $K_{3,3}$-free 2-connected toroidal (non-planar) graphs. Given two disjoint $K_5$-graphs, the graph obtained by identifying an edge of one of the $K_5$'s with an edge of the other is called an {\it $M$-graph} (see Figure~\ref{fig:Mgraph}~(i)), and, when the edge of identification is deleted, an {\it $M^*$-graph} (see Figure~\ref{fig:Mgraph}~(ii)).  

\begin{figure}[h]
	\begin{center}
	\includegraphics*[width=3.6in]{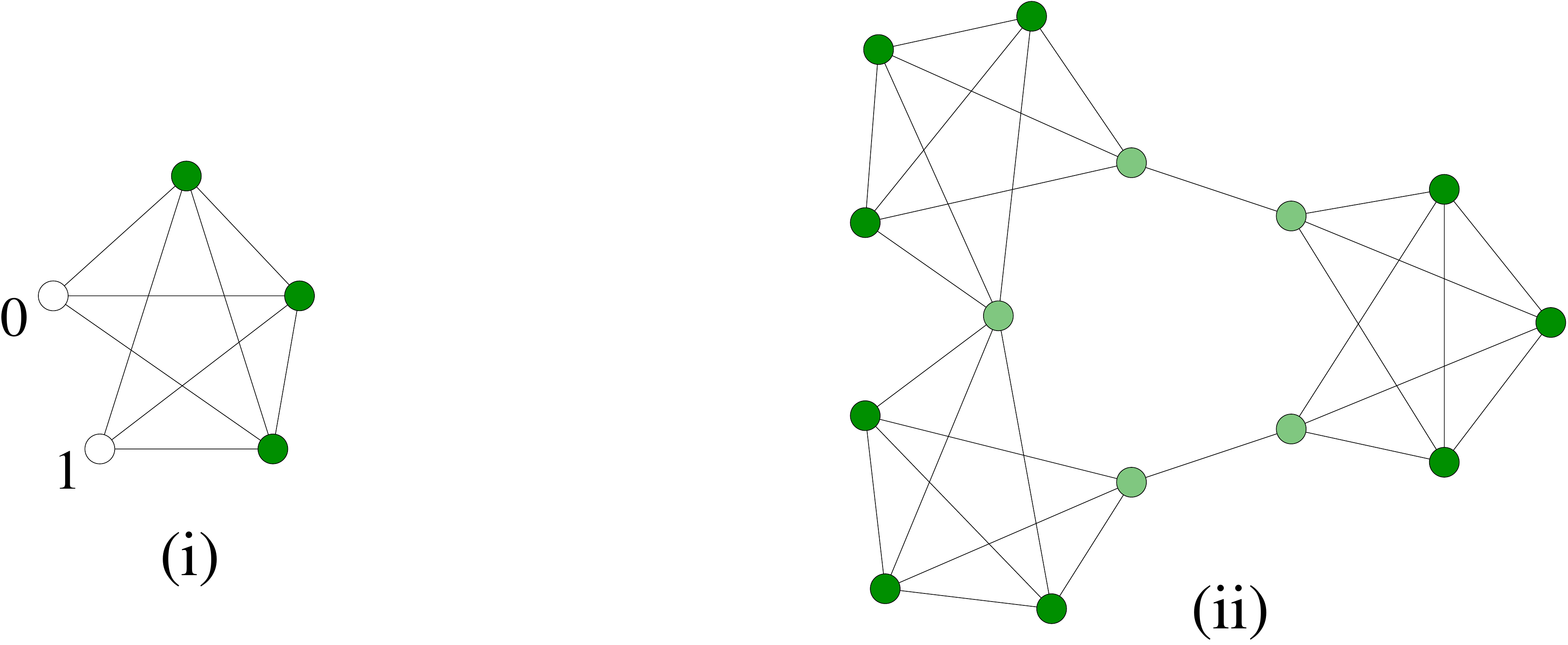}
	\end{center}
\vspace{-5mm}
	\caption{(i) A $(K_{5})_{01}$-network\ \ \ \  (ii) A toroidal crown \label{fig:couronne}}
\end{figure}

A {\it toroidal crown} is a graph $H$ obtained from an unoriented cycle $C_i$, $i\ge 3$, by substituting $(K_{5})_{01}$-networks for some edges of $C_i$ in such a way that no two unsubstituted edges of $C_i$ are adjacent in $H$ (see Figure~\ref{fig:couronne}~(ii) for an example). Denote by $\H$ the class of toroidal crowns.
A {\it toroidal core} is defined as either $K_5$, an $M$-graph, an $M^*$-graph, or a toroidal crown.  
\begin{theorem}[\cite{GLL2}]  The species $\T$ of $K_{3,3}$-free $2$-connected non-planar toroidal graphs can be expressed as a canonical composition
\begin{equation}
\T = \T_C\uparrow \N_P,
\end{equation}
where $\T_C$ denotes the class of toroidal cores, i.e. $\T_C = K_5 + M + M^* + \H$.
\hfill\rule{2mm}{2mm}
\end{theorem}
%

%
In \cite{GLL3} we give explicit formulas for the edge index series 
for the graphs $K_5,\ M,$ and $M^*$ 
and for the species $\H$ of toroidal crowns. 
In order to enumerate unlabelled $K_{3,3}$-free toroidal graphs in $\PP$ or in $\T$ according to the number of vertices and edges, the generating functions $\widetilde{\N_P}(x,y)$ and $\widetilde{\N_P}_\tau(x,y)$ are also required, since we have
\begin{eqnarray}
\widetilde{\PP}(x,y) &=& W_{K_5}[x;\widetilde{\N_P}(x,y); \widetilde{\N_P}_\tau(x,y)], \\
\widetilde{\T}(x,y) &=& W_{\T_C}[x;\widetilde{\N_P}(x,y); \widetilde{\N_P}_\tau(x,y)].
\end{eqnarray}
Using the results of Section~\ref{sec:planar}, we have extended previous tables to 17 vertices for $\widetilde{\PP}(x,y)$ and to 20 vertices for $(\widetilde{\T}-\widetilde{\PP})(x,y)$.  
%
\begin{eqnarray} 
\widetilde{\PP}(x) &=&
{x}^{5}+2\,{x}^{6}+14\,{x}^{7}+102\,{x}^{8}+962\,{x}^{9}+10662\,{x}^{
10}\nonumber\\
&& +139764\,{x}^{11}+2088482\,{x}^{12}+34680722\,{x}^{13}+622943224\,{
x}^{14}\nonumber\\
&& +11854223815\,{x}^{15}+235386309134\,{x}^{16}+4826871283270\,{x
}^{17}  + \cdots \phantom{MMM}
\end{eqnarray} 
\begin{eqnarray} 
(\widetilde{\T}-\widetilde{\PP})(x) &=&
2\,{x}^{8}+11\,{x}^{9}+127\,{x}^{10}+1388\,{x}^{11}+16905\,{x}^{12}\nonumber\\
&& +
214191\,{x}^{13}+2890154\,{x}^{14}+41748279\,{x}^{15}+650024679\,{x}^{
16}\nonumber\\
&& +10888386896\,{x}^{17}+194674234840\,{x}^{18}+3674322404851\,{x}^{
19}\phantom{MMM}\nonumber\\
&& +72412623360105\,{x}^{20}+ \cdots 
\end{eqnarray} 
%
%
%
%
%
\subsection{Homeomorphically irreducible graphs.}
A graph is called \emph{homeomorphically irreducible} if it contains no vertex of degree 2. In order to enumerate these graphs, we can apply the method of Walsh and Robinson (\cite{RobinsonWalsh, Timothyunlabelled}) as follows. 
Any 2-connected graph $G$ is either a series-parallel graph or contains a unique 
2-connected homeomorphically irreducible core $C(G)$, which is different from $K_2$, and unique components $\{N_e\}_{e\in E(C(G))}$ which are series-parallel networks, whose composition gives $G$. 
Let $\B$ be a species of 2-connected graphs. Denote by $I_{\B}$ the class of graphs which are homeomorphically irreducible cores of graphs in $\B$. Also set $\BSP = \B\cap\SPG$ which is the class of series-parallel graphs in $\B$ and let $\DSP$ denote the species of series-parallel networks. 
We then have the following Proposition.
\begin{propos} [\cite{GLL, Timothy}] \label{prop:BBSHBR}
Let $\B$ be a species of 2-connected graphs such that 
\begin{enumerate}
\item
$I_{\B}$ is contained in $\B$,
\item
$\B$ is closed under edge substitution by series-parallel networks.
\end{enumerate}
Then we have 
\begin{equation} \label{eq:BBSHBR}
\B = \BSP + I_{\B}\uparrow \DSP,
\end{equation}
the composition $I_{\B}\uparrow \DSP$ being canonical.
\hfill\rule{2mm}{2mm}
\end{propos}
\begin{propos} [\cite{RobinsonWalsh}] \label{prop:Rinv} 
There exist unique series $\beta(x,y)$ and $\gamma(x,y)$ satisfying
\begin{equation} \label{eq:rhobetagamma}
\rho^+[x,\beta(x,y),\gamma(x,y)] = y, \ \ \ \ \rho^-[x,\beta(x,y),\gamma(x,y)] = y,
\end{equation}
where $\rho^+=W^+_{\DSP}$ and $\rho^-=W^-_{\DSP}$.  Moreover these series are given explicitly by  
\begin{eqnarray} 
\beta(x,y) &=& -1+(1+y)\prod_{j\geq1}(1-x^{2j-1}y^{2j})(1-x^{2j}y^{2j+1})^{-1},
\label{eq:betaexpl}\\
\gamma(x,y) &=& -1+(1+y)\prod_{j\geq1} \frac{(1-x^{4j-3}y^{4j-2})(1+x^{4j-1}y^{4j})}{(1+x^{4j-2}y^{4j-1})(1-x^{4j}y^{4j+1})}. \phantom{MM} \label{eq:gammaexpl}
\end{eqnarray} 
\end{propos}
\proof  We first note that
\begin{eqnarray*} 
\rho^+ &=& b_1 +a_1b_1^2 + a_1b_1^3 + \cdots, \\
\rho^- &=& c_1 + a_1b_2 + a_1b_2c_1 \cdots, \phantom{MM}
\end{eqnarray*} 
where the remaining terms are of higher order in the vertex-cycle variables so that equations (\ref{eq:rhobetagamma}) determine recursively unique series $\beta(x,y)$ and $\gamma(x,y)$. In fact, from (\ref{eq:rhoderho+SP}) and (\ref{eq:rhoderho-SP}), we see that $\beta(x,y)$ and $\gamma(x,y)$ must satisfy
\begin{eqnarray} 
1+y &=& (1+ \beta(x,y))\exp\left(\sum_{i\geq1}  \frac{1}{i} \frac{x^iy^{2i}}{1+x^iy^i}\right),  \\
1+y &=& (1+ \gamma(x,y))\exp\left(\sum_{i\ \mathrm{even}} \frac{1}{i} \frac{x^iy^{2i}}{1+x^iy^i} + \sum_{i\ \mathrm{odd}} \frac{1}{i} \frac{(x^i+x^{2i}y^i)y^{2i}}{1+x^{2i}y^{2i}}\right), \phantom{M}
\end{eqnarray} 
from which (\ref{eq:betaexpl}) and (\ref{eq:gammaexpl}) are readily deduced.
\hfill\rule{2mm}{2mm}
\begin{propos} [\cite{RobinsonWalsh}] \label{prop:WSPGbetagamma} 
For the species $\SPG$ of series parallel graphs, we have
\begin{equation} \label{eq:WSPGbetagamma}
W_{\SPG}[x;\beta(x,y);\gamma(x,y)] = -x^2y^2+ xy(x+xy(1-x))(1-x^4y^4)^{-1}.
\end{equation}
\end{propos}
\proof  Notice that for the edge index series $\sigma^+=W^+_{\S}$ and $\sigma^-=W^-_{\S}$ of series-parallel $s$-networks, we have, using (\ref{eq:sigma}),
\begin{eqnarray} 
\sigma^+[x;\beta(x,y);\gamma(x,y)] &=& \frac{xy^2}{1+xy}, \ \ \ \ \ 
\sigma^-[x;\beta(x,y);\gamma(x,y)]\ =\ \frac{(x+x^2y)y^2}{1+x^2y^2}. \label{eq:sigmabetagamma}\phantom{M}
\end{eqnarray} 
It is then possible to use the dissymmetry formula (\ref{eq:WdissymBisSP}) and the result follows after some simplifications.
\hfill\rule{2mm}{2mm}
\begin{corol} \label{corol:IBtilde}
Let $\B$ be a species of 2-connected graphs such that the hypothesis of Proposition \ref{prop:BBSHBR} are satisfied. Then we have 
\begin{equation} \label{eq:IBtilde}
\widetilde{I_{\B}}(x,y) = (W_\B-W_{\BSP})[x;\beta(x,y);\gamma(x,y)],
\end{equation}
where $\BSP = \B\cap\SPG$.
\end{corol}
%
%
\subsubsection{Example: Planar graphs} 
For the species $\B=\B_P$ of 2-connected planar graphs, we have $\B\cap\SPG=\SPG$.  It follows from Proposition \ref{prop:WSPGbetagamma} and Corollary \ref{corol:IBtilde} that for the species $I_P$ of 2-connected homeomorphically irreducible planar graphs, we have
\begin{equation} \label{eq:IPtildexy}
\widetilde{I_P}(x,y) = W_{\B_P}[x;\beta(x,y);\gamma(x,y)] + x^2y^2 - xy(x+xy(1-x))(1-x^4y^4)^{-1}.
\end{equation}
%

\subsubsection{Example: $K_{3,3}$-free 2-connected graphs}
As seen in Section \ref{sec:K33free}, the species $\B=\B_\F$  associated to the class $\F=\F_P + K_5$ consists of $K_{3,3}$-free 2-connected graphs. Again we have $\SPG\subset\B$ and for the species $I$ of homeomorphically irreducible $K_{3,3}$-free 2-connected graphs we have 
\begin{equation} \label{eq:Itildexy}
\widetilde{I}(x,y) = W_{\B}[x;\beta(x,y);\gamma(x,y)] + x^2y^2 - xy(x+xy(1-x))(1-x^4y^4)^{-1}.
\end{equation}
%
%
%

%
\subsubsection{Example: $K_{3,3}$-free projective planar and toroidal graphs}
For the species $\B=\PP$ and $\B=\T$ of 2-connected $K_{3,3}$-free (non-planar) projective planar and toroidal graphs, respectively, we have 
$\B\cap\SPG=\emptyset$. It follows that for the corresponding species 
$I_{\PP}$ and $I_{\T}$ of homeomorphically irreducible graphs we have
\begin{eqnarray} 
\widetilde{I_{\PP}}(x,y) &=& W_{\PP}[x;\beta(x,y);\gamma(x,y)]
 \nonumber\\
&=& W_{K_5\uparrow\N_P}[x;\beta(x,y);\gamma(x,y)] \nonumber \\
&=& W_{K_5}\left[x;W^+_{\N_P}[x;\beta(x,y);\gamma(x,y)];W^-_{\N_P}[x;\beta(x,y);\gamma(x,y)]\right]. \phantom{M}
\label{eq:IPPtildexy}
\end{eqnarray}
and
\begin{eqnarray} 
\widetilde{I_{\T}}(x,y) &=& W_{\T}[x;\beta(x,y);\gamma(x,y)]
\nonumber \\
&=& W_{\T_C}\left[x;W^+_{\N_P}[x;\beta(x,y);\gamma(x,y)];W^-_{\N_P}[x;\beta(x,y);\gamma(x,y)]\right]. \phantom{M}
\label{eq:ITtildexy}
\end{eqnarray}
%
%
%




\end{document}